\newcommand{\bull}{\vrule height .9ex width .8ex depth -.1ex}
\newcommand{\ppp}{\hfill $\bull$ }
 \documentclass[11pt]{article}
\usepackage{amsfonts}
\usepackage[francais]{babel}
 \author{ Mohammed Larbi Labbi}
   \title{ Courbure riemannienne:\\ 
 variations sur diff\'erentes notions de positivit\'e}

   \date{}
   \usepackage{amsmath}
   \usepackage{amscd}
\newtheorem{theorem}{Th\'eor\`eme}[section]
\newtheorem{corollary}[theorem]{Corollaire}

\newtheorem{proposition}[theorem]{Proposition}
\newtheorem{definition}{D\'efinition}[section]

\begin{document}

\maketitle
\begin{abstract}
Cet article est le texte de  mon habilitation \`a diriger des recherches.
On \'etudie diff\'erentes notions de courbure riemanniennes: 
la $p$-courbure, qui interpole entre courbure scalaire et courbure
 sectionnelle, les courbures de Gauss-Bonnet-Weyl qui 
constituent une autre interpolation allant de la courbure scalaire 
jusqu'\`a l'int\'egrand de Gauss-Bonnet.
 Les $(p,q)$-courbures que nous d\'egageons englobent toutes ces notions.
 On examine ensuite le terme en courbure de la formule classique 
de Weitzenb\"ock. On \'etudie aussi les propri\'et\'es de positivit\'e de 
la $p$-courbure, la seconde courbure de Gauss-Bonnet-Weyl, 
la courbure d'Einstein et de la courbure isotrope.

\vspace{1cm}
 \begin{center}{\bf {Abstract}}\end{center}
This paper is the text of my habilitation thesis. 
We study different notions of Riemannian curvatures: 
The $p$-curvatures interpolate
between the scalar curvature and the sectional curvature, the  
Gauss-Bonnet-Weyl curvatures form another interpolation from
 the scalar curvature to the Gauss-Bonnet integrand. We 
 bring out the
 $(p,q)$-curvatures, which incorporate all the previous curvatures.
We then examine the curvature term which appears in the classical 
Weitzenb\"ock formula.\\
We also study the positivity properties of the $p$-curvatures, the 
second Gauss-Bonnet-Weyl curvature, the Einstein curvature and the 
isotropic curvature.

\end{abstract}

\tableofcontents

\section*{Introduction}

Mes recherches ont port\'e sur la g\'eom\'etrie et l'analyse sur les vari\'et\'es 
riemanniennes compactes de dimension $\geq 4$. Dans ce rapport, 
j'ai regroup\'e mes travaux
suivant quatre th\`emes principaux:
\begin{enumerate}
\item {\bf Les formes doubles:}
Notre objectif initial \'etait l'\'etude du produit ext\'erieur 
sur les formes
doubles, dit produit de "Kulkarni-Nomizu", en vue de pouvoir manipuler des courbures d'ordre 
sup\'erieur. En particulier, les courbures
de Gauss-Bonnet-Weyl, les
$(p,q)$-courbures et les op\'erateurs de courbure de Weitzenb\"{o}ck. Pour
 illustrer cette id\'ee, signalons que l'int\'egrand de Gauss-Bonnet,
qui a la r\'eputation d'\^etre une expression tr\`es compliqu\'ee, s'\'ecrit
 tout simplement  comme une puissance du
tenseur de courbure (vue en tant qu'une forme double). Aussi, l'exemple de l'op\'erateur de courbure
qui apparait dans la formule de Weitzenb\"{o}ck est  significatif, voir chapitre 4.\\
J'ai aussit\^ot constat\'e  que cette \'etude 
 a son int\'er\^et propre. Les formes doubles se pr\'esentent
naturellement comme une g\'en\'eralisation  des formes diff\'erentielles. Il
s'agit alors  de g\'en\'eraliser toute la th\'eorie  classique 
aux formes doubles. C'est un domaine qui \`a mon avis reste encore
 largement inexploit\'e.

\item{\bf Courbures de Gauss-Bonnet-Weyl et g\'en\'eralisations:}
Les courbures de Gauss-Bonnet-Weyl sont les invariants de courbure qui appara\^{\i}ssent
dans la c\'el\`ebre formule des tubes de  Weyl. Elles forment une interpolation entre la 
courbure scalaire usuelle et l'int\'egrand de Gauss-Bonnet.\\
On a d\'emontr\'e que les courbures totales associ\'ees satisfont, comme pour le cas de la courbure 
scalaire totale,
des formules de variation remarquables. Celles-ci donnent naissance \`a des tenseurs 
\`a divergences nulles
du m\^eme  type que le tenseur d'Einstein usuel. Ces r\'esultats
nous ont amen\'e naturellement
\`a des g\'en\'eralisations du probl\`eme de Yamabe pour
 la courbure scalaire et
 des vari\'et\'es
d'Einstein.\\
Par ailleurs, on a g\'en\'eralis\'e la formule d'Avez pour l'int\'egrand de Gauss-Bonnet en dimension
4 aux dimensions sup\'erieures. On en a d\'eduit en particulier, sous certaines hypoth\`eses
g\'eom\'etriques, la non-n\'egativit\'e des courbures de Gauss-Bonnet-Weyl. Ces hypoth\`eses 
seront ensuite caract\'eris\'ees g\'eom\'etriquement \`a l'aide des $(p,q)$-courbures.\\
 Ces derni\`eres g\'en\'eralisent en m\^eme temps 
les $p$-courbures et les courbures de Gauss-Bonnet-Weyl. 
On a d\'emontr\'e \'egalement un th\'eor\`eme 
de Schur
pour ces courbures et  \'etudi\'e leurs propri\'et\'es g\'eom\'etriques.
 
\item{\bf Positivit\'e de la seconde courbure de Gauss-Bonnet-Weyl:}
 La seconde courbure de Gauss-Bonnet-Weyl, not\'ee $h_4$,
 est un invariant scalaire qui est quadratique en $R$: le tenseur de 
courbure de Riemann. Pr\'ecis\'ement, 
$$h_4=\|R\|^2-\|cR\|^2+{1\over 4}\|c^2R\|^2.$$
Une propri\'et\'e importante de la seconde courbure de Gauss-Bonnet-Weyl est qu'elle est
non-n\'egative pour les m\'etriques d'Einstein. D'o\`u l'int\'er\^et majeur de classifier
les vari\'et\'es de dimension $\geq 4$ admettant  une m\'etrique \`a seconde
 courbure de Gauss-Bonnet-Weyl
non-n\'egative. Ceci est un grand projet de recherche et se pr\'esente comme la suite naturelle
du probl\`eme similaire sur la courbure scalaire.\\
On a alors utilis\'e les techniques de chirurgies, submersions  riemanniennes et des actions
de groupe
de Lie  pour construire
des classes de m\'etriques  \`a seconde courbure de Gauss-Bonnet-Weyl
positive. Aussi nous avons pu \'etablir des liens avec la positivit\'e de la courbure 
sectionnelle,
de la courbure isotrope et de la $p$-courbure.
Les r\'esultats obtenus montrent une forte analogie entre le comportement de
la positivit\'e de la courbure scalaire et celle de $h_4$. Voici quelques exemples:
\begin{enumerate}
\item La classe des vari\'et\'es compactes \`a $h_4>0$ est stable sous chirurgies en
codimension $>4$. 
\item Pour les vari\'et\'es compactes qui sont les espaces totaux pour des
submersions riemanniennes, il suffit que les fibres soient \`a $h_4>0$ pour
pouvoir construire une m\'etrique \`a $h_4>0$ sur l'espace total.
\item Soit $n>5$. Tout groupe de pr\'esentation finie peut \^etre r\'ealis\'e comme le
groupe fondamental d'une vari\'et\'e de dimension $n$ \`a $h_4>0$.
\end{enumerate}
En contre partie, il y a des diff\'erences  subtiles entre la positivit\'e de 
ces deux courbures. En effet,
on sait qu'une vari\'et\'e compacte
admettant une m\'etrique \`a courbure sectionnelle n\'egative ne peut pas porter de 
m\'etriques \`a courbure scalaire positive (c'est un r\'esultat de Gromov). En revanche,
il n'y a pas de probl\`eme d'incompatibilit\'e entre la n\'egativit\'e de la 
courbure sectionnelle 
et la positivit\'e de $h_4$. Ceci est illustr\'e tout simplement par
 l'exemple d'une
vari\'et\'e \`a courbure constante n\'egative arbitraire!.

\item{\bf Courbure isotrope positive:}
La courbure isotrope a \'et\'e introduite par Micallef et Moore  pour les 
vari\'et\'es
de dimension $\geq 4$. Cette courbure joue le m\^eme r\^ole dans l'\'etude de
la stabilit\'e des 2-sph\`eres
harmoniques et dans celle des surfaces minimales, que celui exerc\'ee par
 la courbure sectionnelle
dans l'\'etude de la stabilit\'e des g\'eod\'esiques.\\
L'\'existence d'une m\'etrique riemannienne \`a courbure isotrope positive
sur une vari\'et\'e compacte entra\^{\i}ne l'annulation des groupes d'homotopies
$\pi_i(M)$ pour $2\leq i\leq [n/2]$.\\
On a \'etabli une relation \'etroite entre la courbure isotrope et la $(n-4)$-courbure,
 o\`u $n$ est la dimension
de la vari\'et\'e en question.
En effet, ces deux courbures  se distinguent seulement par un terme en la courbure de Weyl,
 et par cons\'equent, 
elles co\"{\i}ncident dans le cas conform\'ement plat. J'ai alors utilis\'e cette analogie pour
adapter mes r\'esultats de th\`ese sur la $p$-courbure au cas de la courbure isotrope.
 \end{enumerate}
\newpage
\part{Formes doubles et courbures riemanniennes}
\section{Les  formes doubles}
Le tenseur covariant  de courbure de Riemann est antisym\'etrique par
rapport aux deux premiers variables ainsi que par rapport aux  deux derniers.
On pourra alors le consid\'erer comme une forme bilin\'eaire sur l'espace des
 bivecteurs. C'est un exemple typique d'une forme double sym\'etrique d'ordre (2,2).\\
Un deuxi\`eme exemple est le tenseur de  Thorpe \cite{Thorpe,Stehney}. C'est 
 un tenseur de courbure
d'ordre sup\'erieur ayant pour courbure sectionnelle les courbures
de Gauss-Kronecker.  Celui-ci  pr\'esente les m\^emes propri\'et\'es
d'antisym\'etrie par rapport aux $p$-premiers  variables et par rapport aux
 $p$-derniers. C'est  une forme double sym\'etrique d'ordre $(p,p)$.
  Un troisi\`eme
exemple significatif  
de forme double sym\'etrique  est l'op\'erateur de 
courbure de  Weitzenb\"{o}ck. Ce dernier sera \'etudi\'e dans le chapitre 4.  \\
D'autres exemples naturels de formes doubles en g\'eom\'etrie riemannienne 
sont la courbure de Ricci,
 le tenseur
 d'Einstein, le tenseur de courbure de Weyl, le tenseur de Schouten,
  la seconde
  forme fondamentale,.... \\
Le produit ext\'erieur usuel sur les formes
 se g\'en\'eralise naturellement aux
formes doubles. On obtient alors le produit dit de "Kulkarni-Nomizu''.\\
Ce produit a le m\'erite de rendre maniable des
expressions naturelles, mais compliqu\'ees,  en la courbure
de Riemann. Un  exemple signifiant est  l'int\'egrand de Gauss-Bonnet.
Il  s'\'ecrit tout
simplement, \`a une constante pr\`es, comme une puissance
 du tenseur de courbure de Riemman. Un autre exemple est l'op\'erateur
 de courbure de la formule de Weitzenb\"{o}ck, voir chapitre 4.\\
Les  formes doubles ont \'et\'e  introduites pour la premi\`ere fois par De Rham,
 ensuite reprises et d\'evolopp\'ees dans les ann\'ees 70 par
Kulkarni \cite{Kulkarni},
Nomizu \cite{Nomizu}, Gray\cite{Gray}, Kowalski \cite{Kowalski}.... Elles sont
 aussi \'etudi\'ees
 en physique math\'ematique, voir les travaux r\'ecents  de Senovilla \cite{Senovilla} et
 les 
r\'ef\'erences qui y sont cit\'ees.\\
Dans ce chapitre, on se propose de  pr\'esenter notre contribution dans ce sujet. Il s'agit
essentiellement
de r\'esultats techniques mais d'utilit\'e majeure pour la suite de notre \'etude.   \\
Pour fixer les notations commen\c cons par les d\'efinitions:\\
Soit $(V,g)$ un espace vectoriel Euclidien (orient\'e) r\'eel de dimension n.
Dans la suite on identifiera, sans le signaler,
 les espaces vectoriels et leurs
duaux via les structures  Euclidiennes. Soit
  $\Lambda^{*}V=\bigoplus_{p\geq 0}\Lambda^{*p}V$ (resp.
   $\Lambda V=\bigoplus_{p\geq 0}\Lambda^{p}V$)  l'alg\`ebre
 des $p$-formes (resp. $p$-vectors) sur $V$. En prenant des produits
 tensoriels,
  on d\'efinit l'espace des formes doubles
 ${\cal D}= \Lambda^{*}V\otimes \Lambda^{*}V=\bigoplus_{p,q\geq 0}
  {\cal D}^{p,q}$ o\`u $  {\cal D}^{p,q}= \Lambda^{*p}V \otimes
   \Lambda^{*q}V$.  C'est une alg\`ebre  bi-gradu\'ee associative,
    o\`u la multiplication est not\'ee par un point, celui ci sera omis
     d\`es que convenable.
\par\noindent
    Rappelons que tout \'el\'ement du produit tensoriel
    $  {\cal D}^{p,q}= \Lambda^{*p}V \otimes \Lambda^{*q}V$
    peut \^etre canoniquement identifi\'e \`a une forme  bilin\'eaire
 $\Lambda^pV\times\Lambda^qV\rightarrow {\bf R}$. Ceci \'etant une forme
 multilin\'eaire
 antisym\'etrique par rapport aux $p$-premiers arguments ainsi qu'aux
  $q$-derniers.
    Sous cette identification, le produit de $\omega_1\in {\cal D}^{p,q}$ et
$\omega_2\in {\cal D}^{r,s}$
est donn\'e par
\begin{equation*}
\begin{split}
\label{eps:prod}
 &\omega_1.\omega_2(x_1\wedge...\wedge x_{p+r},y_1\wedge...\wedge y_{q+s})\\
&= {1\over p!r!s!q!}\sum_{\sigma\in S_{p+r}, \rho\in S_{q+s}}
\epsilon(\sigma)\epsilon(\rho)
\omega_1(x_{\sigma(1)}\wedge...\wedge x_{\sigma(p)};y_{\rho(1)}
\wedge...\wedge y_{\rho(q)})\\
&\phantom{...mmmmmmmmmmm}
\omega_2(x_{\sigma(p+1)}\wedge...\wedge x_{\sigma(p+r)};y_{\rho(q+1)}
\wedge...\wedge y_{\rho(q+s)}).
\end{split}
\end{equation*}
En particulier, le produit du produit scalaire $g$
avec lui m\^eme $k$-fois d\'etermine
le produit scalaire canonique sur $\Lambda^{p}V$. Pr\'ecis\'ement on a
$$g^k(x_1 \wedge...\wedge x_k,y_1\wedge...\wedge y_k)=k!\det[g(x_i,y_j)].$$
\subsection{L'application de multiplication par $g^k$}
L'application de  multiplication par $g^k$ dans ${\cal D}$ joue un r\^ole
 fondamental
dans l'\'etude des formes doubles.
On a montr\'e le r\'esultat suivant qui g\'en\'eralise un lemme d\^u \`a Kulkarni
\cite{Kulkarni} dans le cas
$k=1$:
\begin{proposition}[\cite{Labbi7}]\label{mult:injec} L'application de multiplication
 par $g^k$ est
 une application
injective sur $D^{p,q}$
d\`es que $p+q+k<n+1$.\end{proposition}
Ce r\'esultat  nous a permis d'une part de simplifier des calculs autrement
tr\`es compliqu\'es, voir chapitre 3 et d'autres exemples de calculs
simplifi\'es dans \cite{Labbi7, Labbi11},
et d'autre part de mieux comprendre la structure de l'alg\`ebre des formes
 doubles. \\

Une deuxi\`eme application fondamentale dans ${\cal D}$ est l'application
de contraction $c$. Elle envoie ${\cal D}^{p,q}$ sur ${\cal D}^{p-1,q-1}$ et
elle est l'adjointe  de la multiplication par $g$ comme on le verra
 dans la section
 ci-dessous.\\
 Ces deux applications ne commutent pas en g\'en\'eral. Dans \cite{Labbi7}
 on a \'etabli  apr\`es un calcul d\'elicat, une formule explicite
  pour la diff\'erence $c^kg^l-g^lc^k$. Cette formule a \'et\'e l'ingr\'edient
  principal
 dans l'\'elaboration de la proposition pr\'ec\'edente.\\
 Signalons au passage que,  contrairement \`a toute impression spontan\'ee, 
notre  proposition ne peut pas \^etre obtenue
  par de
 simples applications successives du lemme de Kulkarni.\\
Passons maintenant au produit scalaire sur ${\cal D}$.
Le produit scalaire canonique sur $\Lambda^{*p}V$
  induit naturellement un produit scalaire sur les espaces
  $D^{p,q}=\Lambda^{*p}V \otimes \Lambda^{*q}V$. On le notera  par   $ <,>$.
  \par
On \'etend $<,>$ \`a ${\cal D}$ en posant  $D^{p,q} \perp D^{r,s}$
si $p\not = r$ ou si $q\not =s$.
\\
Dans la proposition suivante on a montr\'e que la contraction dans ${\cal D}$ est
l'adjointe  de
la multiplication par $g$:\\
\begin{proposition}[\cite{Labbi7}]\label{theo:gc} Pour tout $\omega_1, \omega_2\in {\cal D}$  on a
\begin{equation}
\label{adj:gc}
<g\omega_1,\omega_2>=<\omega_1,c\omega_2>.
\end{equation}
En particlier, pour tout $k\geq 1$, on a
$<g^k\omega_1,\omega_2>=<\omega_1,c^k\omega_2>$.
\end{proposition}
\subsection{L'op\'erateur de Hodge g\'en\'eralis\'e}
On suppose que l'espace $V$ est orient\'e. L'op\'erateur de Hodge $*:\Lambda^{p}V^*\rightarrow \Lambda^{n-p}V^*$
 s'\'etend de mani\`ere naturelle pour donner un op\'erateur lin\'eaire  $*:{ D}^{p,q}
 \rightarrow { D}^{n-p,n-q}$. Pour $\omega=\theta_1\otimes
\theta_2$ on pose
$$*\omega=*\theta_1\otimes *\theta_2.$$
Notons que  $*\omega(.,.)=\omega(*.,*.)$ en tant que forme bilin\'eaire.
 Plusieurs propri\'et\'es
de l'op\'erateur de Hodge usuel se g\'en\'eralisent alors
naturellement, voir \cite{Labbi7,Labbi9}.\\
L' op\'erateur de Hodge g\'en\'eralis\'e s'est r\'ev\'el\'e un outil tr\`es utile dans notre 
\'etude et
 nous a permis, en particulier, d' \'etablir
une deuxi\`eme formule simple reliant
l'application de contraction $c$ \`a la multiplication par $g$:\\
\begin{proposition}[\cite{Labbi7}] Pour tout $\omega \in  {\cal D}$, on a
\begin{equation}
\label{gstar:c}
g\omega=*c*\omega.
\end{equation}
 En particulier, pour tout $k\geq 1$ 
et $\omega \in D^{p,p}$, on a
$g^k\omega=*c^k*\omega$.
\end{proposition}
Mentionnons au passage que l'op\'erateur de Hodge peut \^etre \'etendu aux formes doubles
selon deux autres  mani\`eres suppl\'ementaires.
 Pour $\omega$ comme ci-dessus on pose 
$*_1\omega=*\theta_1\otimes \theta_2$ et $*_2\omega=\theta_1\otimes *\theta_2.$ 
Les op\'erateurs ainsi obtenus ont l'inconv\'enient de ne pas conserver la symm\'etrie mais 
ils ont leurs int\'er\^ets propres. Voir les travaux r\'ecents de \cite{Senovilla} 
pour des applications en physique math\'ematique
de ces diff\'erentes extensions.

\subsection{D\'ecomposition orthogonale}
Notre objectif ici est de g\'en\'eraliser, aux tenseurs de courbures d'ordre sup\'erieur,
la  d\'ecomposition standard du
tenseur de courbure de Riemann $R$, qui rappelons le affirme que:
$$R=W+{1\over n-2}(c(R)-{1\over n}g.c^2(R))g+{1\over 2n(n-1)}c^2(R).g^2.$$
O\`u  $W$, $cR$ et $c^2R$ d\'esignent r\'espectivement les courbures de Weyl,  Ricci et
la courbure scalaire. 
\par\noindent
On se place, comme on l'a fait jusqu'ici, dans le cadre alg\'ebrique.
Remarquons  que pour  $\omega_1\in \ker c$, $g\omega_2\in g
D^{p,q}$ et en utilisant  (\ref{adj:gc}), on a
$<\omega_1,g\omega_2>=<c\omega_1,\omega_2>=0.$
On obtient alors (pour des raisons de dimension) la d\'ecomposition orthogonale 
$D^{p+1,q+1}={\rm Ker}\, c\oplus gD^{p,q},$
et ceci pour tout  $p,q\geq 0$ tels que $p+q\leq n-1$. Si on note  par $E^{p,q}$ le noyau
 $\ker c\subset D^{p,q}$, on obtient
par cons\'equent  la d\'ecomposition orthogonale suivante
de $D^{p,q}$:
\begin{equation}\label{ort:decom}
D^{p,q}=E^{p,q}\oplus gE^{p-1,q-1}\oplus g^2E^{p-2,q-2}\oplus ...
\oplus g^rE^{p-r,q-r},
\end{equation}
o\`u $r={\min \{p,q\}}$.
Pour une forme double donn\'ee, on a pu  \'etablir
 une formule explicite pour ses diff\'erentes composantes
suivant la  d\'ecomposition pr\'ec\'edente comme suit:
\begin{theorem}[\cite{Labbi7}]\label{omk:form}
Suivant la d\'ecomposition orthogonale (\ref{ort:decom}),
 toute  $\omega\in D^{p,p}$ se d\'ecompose comme suit
$$\omega=\omega_p+g.\omega_{p-1}+...+g^p.\omega_0,$$
o\`u $\omega_0={(n-p)!\over p!n!}c^p(\omega)$, et pour $1\leq k\leq p$ on a
\begin{equation*}
\omega_k={(n-p-k)!\over (p-k)!(n-2k)!}\left[ c^{p-k}(\omega)+
\sum_{r=1}^k{(-1)^r\over {\prod_{i=0}^{r-1}(n-2k+2+i)}}
{g^r\over r!}c^{p-k+r}(\omega)\right].
\end{equation*}
\end{theorem}
On retrouve, en particulier, la  d\'ecomposition standard cit\'ee ci-dessus du
tenseur de courbure de Riemann pour $p=2$ et $\omega=R$.\\

Notons par ailleurs  que pour $\omega\in D^{p,p}$, la composante $\omega_p={\rm con}\, \omega$
 d\'epend
seulement de la classe conforme de la m\'etrique $g$.
 Elle g\'en\'eralise alors le tenseur de courbure de
Weyl. Cette constatation  a amen\'e Nasu \`a \'etudier dans \cite{Nasu2}
une g\'en\'eralisation des vari\'et\'es
conform\'ement plates: les vari\'et\'es $q$-conform\'ement plates. Celles-ci sont caract\'eris\'ees
par l'annulation de la composante $\omega_{2q}$ du tenseur de courbure de Gauss-Kronecker $R^q$.\\
Enfin, signalons  que la d\'ecomposition pr\'ec\'edente n'est pas en g\'en\'eral
irr\'eductible sous l'action
du groupe orthogonal. Ceci \'etant essentiellement d\^u au fait que les sous-espaces
v\'erifiant la premi\`ere
identit\'e de Bianchi sont invariants, voir \cite{Kulkarni}.

\subsection{L'alg\`ebre des structures de courbure}
Notons par ${\cal C}^p$ les \'el\'ements sym\'etriques de  ${D}^{p,p}$. Suivant  Kulkarni, on
 d\'efinit l'alg\`ebre des structures de courbure comme la sous-alg\`ebre  commutative
 ${\cal C}=\bigoplus_{p\geq 0}{\cal C}^p$.
 \par\medskip
Rappelons que la premi\`ere somme de Bianchi,
 not\'ee ici par ${\cal B}$, envoie  ${\cal D}^{p,q}$ dans
 ${\cal D}^{p+1,q-1}$. Le noyau ${\rm ker}{\cal B}$ est stable par la
multiplication dans ${\cal D}$. On appelera alors la sous-alg\`ebre commutative
${\cal C}_1={\cal C}\cap {\rm ker}{\cal B}$, l'alg\`ebre des structures de courbure v\'erifiant
 la premi\`ere identit\'e de Bianchi. \\
Notons par $E_1^p$ les \'el\'ements \`a trace nulle dans ${\cal C}_1^p$.
Il est utile de mentionner ici que Kulkarni \cite{Kulkarni}
a d\'emontr\'e que $E_1^p$ est irr\'eductible sous l'action du groupe orhogonal.
Il en a par cons\'equent  \'etabli la  d\'ecomposition orthogonale de ${\cal C}_1^p$
 en composantes irr\'eductibles comme suit:
$${\cal C}_1^p=E_1^p\oplus gE_1^{p-1}\oplus g^2E_1^{p-2}\oplus ...
\oplus g^pE_1^0.$$
Comme pour le cas du  tenseur de courbure de Riemann,
toute  structure de courbure $\omega\in {\cal C}_1$ est  compl\`etement d\'etermin\'ee
 \`a partir
de sa courbure sectionnelle. De plus, une structure de courbure $\omega\in {\cal C}^p$ 
satisfait 
la premi\`ere identit\'e de Bianchi si et seulement s'il en est de m\^eme pour $*\omega$.
Nous avons  utilis\'e ces constatations pour \'etablir une formule explicite pour
 l'op\'erateur de Hodge. Le r\'esultat est:

\begin{theorem}[\cite{Labbi7}]
Soit $\omega\in {\cal C}_1^p$ et $1\leq p\leq n$, alors on a
\begin{equation*}
*\omega=\sum_{r=\max\{0,2p-n\}}^p
{(-1)^{r+p}\over r!}{g^{n-2p+r}\over (n-2p+r)!}c^r\omega.
\end{equation*}
De plus, suivant la  d\'ecomposition (\ref{ort:decom}),
 pour $\omega=\sum_{i=0}^p g^{p-i}\omega_i$ on a:
\begin{equation*}\label{star:omi}
*\omega=\sum_{i=0}^{\min\{p,n-p\}}(p-i)!(-1)^i
{1\over (n-p-i)!}g^{n-p-i}\omega_i.
\end{equation*}
\end{theorem}

Remarquons que pour $n=2p$, l'op\'erateur de Hodge g\'en\'eralis\'e ne d\'epend que de la
 classe conforme
de la m\'etrique $g$. En effet, les termes $g^rc^r\omega$ restent inchang\'es apr\`es un changement conforme
de la m\'etrique.\\

\subsection{La seconde somme de Bianchi et l'op\'erateur Hessien g\'en\'eralis\'e}
Soit $(M,g)$ une vari\'et\'e riemannienne  de dimension  $n$ et $T_mM$ son espace tangent au
 point $m\in M$. On notera aussi par  $D^{p,q},{\cal C}^p,{\cal C}_1^p...$
les fibr\'es vectoriels sur $M$ ayant pour fibres en $m$, les espaces
$D^{p,q}(T_mM)$, ${\cal C}^p(T_mM)$, ${\cal C}_1^p(T_mM) ...$.\\
On remarque que tous les r\'esultats alg\'ebriques pr\'ec\'edents sont encore valables
sur  l'anneau des sections de ces fibr\'es. Le produit scalaire des sections \'etant bien
 \'evidemment
 le produit
scalaire int\'egral.\par\medskip\noindent
La  seconde somme de  Bianchi  $D$ envoie $ D^{p,q}$ dans  $ D^{p+1,q}$. Sa
restriction  \`a ${\cal D}^{p,0}$ co\"{\i}ncide avec $-d$, o\`u $d$ est
 l'op\'erateur de d\'erivation ext\'erieur sur les  $p$-formes.
Une deuxi\`eme extension naturelle de  $d$ est l'op\'erateur $\tilde D$.
 Il est analogue \`a $D$ mais en revanche  envoie
$ D^{p,q}$ sur $ D^{p,q+1}$. Pr\'ecis\'ement, pour  $\omega\in{\cal D}^{p,q}$, on pose
$$(\tilde{D}\omega)(x_1\wedge...\wedge x_{p},y_1\wedge...\wedge y_{q+1})=
\sum_{j=1}^{q+1}(-1)^j{\nabla_{y_j}\omega}(x_1\wedge...\wedge x_{p},
 y_1\wedge... \wedge \hat{y_j}\wedge ...\wedge y_{q+1}).$$
 Notons qu'en g\'en\'eral, contrairement \`a  la d\'erivation ext\'erieure ordinaire,
 on n'a ni  $D^2=0$ ni ${\tilde D}^2=0$.\\
On introduit maintenant l'op\'erateur  $\delta=c\tilde D+\tilde D c$ qui est 
d\'efini sur les formes doubles et
g\'en\'eralise l'op\'erateur  $\delta$ classique  sur les formes.
Par un calcul direct on a  montr\'e que:
\begin{proposition}[\cite{Labbi9}]
Soient $(M,g)$ une vari\'et\'e riemannienne  orient\'ee, $D$ sa seconde somme de Bianchi, et
  $ *$  l'op\'erateur g\'en\'eralis\'e de Hodge. Alors,
pour toute $(p,q)$-forme $\omega$ sur $M$ telle que  $p\geq 1$ on a
\begin{equation*}\label{stard:star}\delta \omega  =*D* \omega.
 \end{equation*}
De plus, si $M$ est compacte, l'op\'erateur $(-1)^{p+1+(p+q)(n-p-q)}\delta:D^{p+1,q}\rightarrow D^{p,q}$ est
l'adjoint formel de  l'op\'erateur $D$.

\end{proposition}
\par\medskip\noindent
Cette  derni\`ere proposition \'etend aux formes doubles  un r\'esultat classique pour les
 formes diff\'erentielles.\\
D'une mani\`ere similaire, on  introduit l'op\'erateur $\tilde \delta=cD+Dc$. Alors pour toute
$(p,q)$-forme $\omega$ avec $q\geq 1$ on montre que
$\tilde \delta \omega=*\tilde D *\omega.$
Aussi, pour une vari\'et\'e compacte, l'adjoint formel de  $\tilde D$ est l'op\'erateur
$(-1)^{q+1+(p+q)(n-p-q)}\tilde\delta:D^{p,q+1}\rightarrow D^{p,q}.$
\par\medskip\noindent
Finalement, en combinant les op\'erateurs pr\'ec\'edents, on d\'efinit
l'op\'erateur $D\tilde D+\tilde D D$ qui envoie, pour tout $p$,  
${\cal C}^p$ dans ${\cal C}^{p+1}$. Il g\'en\'eralise l'op\'erateur Hessien
usuel sur les fonctions ($p=0$).\\
Il est remarquable que  ce Hessien g\'en\'eralis\'e appara\^{\i}t naturellement
 dans la premi\`ere variation
du tenseur de courbure de Riemann une fois consid\'er\'e comme une forme double sym\'etrique,
voir chapitre 2. Son adjoint formel  est donn\'e, d'apr\`es ce qui pr\'ec\`ede, par
\begin{equation}\label{adjoint:tildedd}
(-1)^{n-p-q-1}[\tilde\delta \delta+\delta\tilde\delta]  =(-1)^{n(p+q)}
*[\tilde DD+D\tilde D]*:
 {\cal C}^{p+1}\rightarrow {\cal C}^{p}.
 \end{equation}

\section{Les courbures de Gauss-Bonnet-Weyl}
\subsection{Introduction}
Les fonctions sym\'etriques paires sur les valeurs propres de la seconde forme fondamentale
d'une hypersurface de l'espace Euclidien sont intrins\`eques. Elles
 peuvent alors \^etre d\'efinies
pour n'importe quelle vari\'et\'e riemannienne et on les  appelle
 alors courbures de Gauss-Bonnet-Weyl. Elles forment une interpolation entre la courbure
scalaire usuelle et l'int\'egrand de Gauss-Bonnet et  appara\^{\i}ssent toutes
dans la c\'el\`ebre formule des tubes de Weyl \cite{Weyl}.\\
 Ces courbures  g\'en\'eralisent naturellement la courbure scalaire. De  plus, elles ont
le m\'erite de satisfaire une formule de variation semblable \`a celle
satisfaite par la courbure scalaire. Cette formule donne naissance \`a des
tenseurs d'Einstein g\'en\'eralis\'es, dits d'Einstein-Lovelock. Ils sont \'egalement,
 comme le tenseur d'Einstein usuel, \`a divergence nulle.\\
On va tout d'abord commencer par la d\'efinition pr\'ecise des ces courbures, ensuite on 
 donnera
quelques exemples pour mieux situer ces invariants parmi les autres courbures bien connues.\\
Soit $R$ le tenseur de courbure de Riemann d'une vari\'et\'e riemannienne
$(M,g)$ de dimension $n$. On note par $R^q$ le produit de $R$ avec lui m\^eme $q$-fois dans
l'anneau de courbure. Soulignons que
les tenseurs $R^q$ co\"{\i}ncident avec
les tenseurs de courbure de Gauss-Kronecker. \\
Il n'est pas inutile de rappeler  que  $R^q$, ainsi que tous les tenseurs  $g^pR^q$ et $*g^pR^q$,
satisfont la premi\`ere identit\'e de Bianchi.
\par\medskip\noindent
{\sl Pour $1\leq 2q \leq n$, la $(2q)$-i\`eme coubure de Gauss-Bonnet-Weyl,
 not\'ee par $h_{2q}$, est la
 contraction compl\`ete du tenseur  $R^q$. Pr\'ecis\'ement, on pose
  $$h_{2q}={1\over (2q)!}c^{2q}R^q.$$
  }\par\medskip\noindent
La d\'enomination de ces courbures est justifi\'ee par les faits suivants:\\
Notons que  $h_2={1\over 2}c^2R$ est la moiti\'e de la courbure scalaire.
Pour
 $n$ paire $(n=\dim M)$, la courbure   $h_n$ co\"{\i}ncide,
  \`a un facteur constant pr\`es,
 avec l'int\'egrand de  Gauss-Bonnet. De plus, pour tout $k$,
 $1\leq 2k \leq n$, la courbure  $h_{2k}$ co\"{\i}ncide, \`a un facteur constant  pr\`es,
 avec les invariants de courbures qui appara\^{\i}ssent naturellement dans la formule
 des tubes de Weyl. \\
 Signalons enfin qu'en utilisant la formule explicite de l'op\'erateur de Hodge
 g\'en\'eralis\'e, on peut r\'e\'ecrire ces courbures comme suit:
 $$h_{2q}= *{1\over (n-2q)!}g^{n-2q}R^q.$$
\par\medskip\noindent
Voici, ci-dessous quelques exemples, pour plus de d\'etailles voir \cite{Labbi7,Labbi8}:
\par\medskip\noindent
\begin{itemize}
\item Soit  $(M,g)$ une vari\'et\'e riemannienne \`a courbure sectionnelle
constante $\lambda$, alors pour tout $q$, la courbure
 $h_{2q}$ est constante et \'egale \`a
$$h_{2q}={\lambda^qn!\over 2^q(n-2q)!}.$$
\item Soient $(M,g)$  une  hypersurface de l'espace  Euclidien et
 $\lambda_1\leq \lambda_2\leq ...\leq \lambda_n$
 les valeurs propres de sa seconde forme fondamentale $B$
en un point donn\'e.
Alors en ce point on a,
 $$h_{2q}=
 {(2q)!\over 2^q}\sum_{1\leq i_1<...<i_{2q} \leq n}\lambda_{i_1}
 ...\lambda_{i_{2q}}.$$
Par cons\'equent, les courbures $h_{2q}$ co\"{\i}ncident, \`a un facteur constant pr\`es,
avec les fonctions sym\'etriques $s_{2q}$ des valeurs
propres de  $B$.

 \item  Soit $(M,g)$  une vari\'et\'e riemannienne conform\'ement plate.
Il est bien connu que son tenseur de courbure est alors  d\'etermin\'e par le tenseur de
Schouten $h$, pr\'ecis\'ement $R=gh$.\\
 Un calcul direct montre que les courbures $h_{2q}$ co\"{\i}ncident,
\`a un facteur constant pr\`es, avec les fonctions sym\'etriques des valeurs
propres de  $h$. Ces courbures ont \'et\'e r\'ecemment \'etudi\'ees par Gursky,
 Viaclovsky, Trudinger,...,
voir \cite{Gursky,Trudinger}. Ils ont, en particulier,  prouv\'e une g\'en\'eralisation du
probl\`eme de Yamabe  pour ces courbures.
\end{itemize}
 \par\medskip
\subsection{Les tenseurs d'Einstein-Lovelock}
Apr\`es la m\'etrique elle m\^eme, le tenseur d'Einstein est,
\`a un facteur constant  pr\`es, l'unique combinaison lin\'eaire de la m\'etrique
et de sa courbure de Ricci  \`a \^etre  \`a
 divergence nulle.\\
En rempla\c cant la courbure de Ricci $cR$ par la courbure de Ricci
 g\'en\'eralis\'ee $c^{2q-1}R^q$ et en imposant les m\^emes conditions cit\'ees ci-dessus on obtient
 le tenseur d'Einstein-Lovelock d'ordre $2q$.\\
Ces tenseurs peuvent \^etre introduits d'une autre mani\`ere, comme suit:\\
Rappelons d'abord qu'on obtient le tenseur d'Einstein,
ou plus pr\'ecis\'ement sa courbure sectionnelle dans une direction $v$,
 en contractant compl\`etement
le tenseur de courbure de Riemann dans les directions orthogonales \`a $v$.
 D'une mani\`ere similaire,  en rempla\c cant
le tenseur de
Riemann par le tenseur de Gauss-Kronecker on obtient les tenseurs d'Einstein-Lovelock.\\
La d\'efinition pr\'ecise est comme suit:
\begin{definition}
Le tenseur  d'Einstein-Lovelock d'ordre $2q$, not\'e $T_{2q}$,
 est d\'efini par
\begin{equation}\label{t2q:def}
T_{2q}=h_{2q}g-{1\over (2q-1)!}c^{2q-1}R^q.
\end{equation}
\end{definition}
Le tenseur ${1\over (2q-1)!}c^{2q-1}R^q$ est l'analogue du tenseur de Ricci
pour le tenseur $R^q$. Remarquons alors l'analogie avec
le tenseur d'Einstein. En particulier,
pour $q=1$, on retrouve ce dernier i.e. 
 $T_2={1\over 2}c^2Rg-cR$.\\
 Notons que ces tenseurs ne satisfont pas  g\'en\'eralement la deuxi\`eme identit\'e
 de Bianchi, en revanche ils sont \`a divergences nulles.\\
La nullit\'e de la divergence pour les tenseurs d'Einstein-Lovelock 
est une cons\'equence du fait g\'en\'eral 
suivant: si une fonctionnelle riemannienne admet un gradient $L^2$, ce gradient est
\`a divergence nulle. Cette remarque est d\^ue \`a D. Bleecker
et remonte en fait \`a Hilbert, voir le chapitre 4 de Besse  \cite{Besse}.

\subsection{Une formule de variation}
 Soit  $M$ une vari\'et\'e $C^\infty$ (orient\'ee) et compacte  de  dimension $n$, et soit
  ${\cal M}$ l'espace des
 m\'etriques riemanniennes $C^\infty$  sur $M$  muni d'une   $L^2$-norme naturelle de Sobolev.
 Celui-ci nous permettra de parler de fonctionnelles diff\'erentiables
 ${\cal M}\rightarrow {\bf R}$.
Une fonctionnelle  $F:{\cal M}\rightarrow {\bf R}$ est dite 
riemannienne si elle est
  invariante sous l'action du groupe des  diff\'eomorphismes.
 On dit que $F$
 a un  gradient en $g$ s'il existe un tenseur sym\'etrique  
$a\in {\cal C}^1$
  tel que
 pour tout tenseur sym\'etrique $h\in {\cal C}^1$ on a
 $$F_g'h={d\over dt}\mid_{t=0}F(g+th)=<a,h>,$$
o\`u ${\cal C}^1$ repr\'esente l'espace des tenseurs sym\'etriques dans
$\Lambda^{*}M \otimes\Lambda^{*}M$ et
$<,>$  le produit scalaire int\'egral.
\par\medskip\noindent
Rappelons que la fonctionnelle classique de courbure scalaire 
totale est d\'efinie par
 $$S(g)=\int_M c^2R \mu_g,$$
 o\`u $c^2R$ repr\'esente la courbure scalaire de $g$
  et $\mu_g$  l'\'el\'ement de
 volume de $g$.  Les points critiques de cette fonctionnelle, une fois restreinte \`a
  ${\cal M}_1=\{ g\in {\cal M}:{\rm vol}(g)=1\}$,
 sont les m\'etriques d'Einstein. \\
La fonctionnelle riemannienne suivante g\'en\'eralise
 d'une mani\`ere naturelle $S$:
$$H_{2k}(g)=\int_Mh_{2k}\mu_g,$$

Remarquons que pour $k=1$, $H_2=S/2$ est la moiti\'e  de $S$.
  Aussi, si la dimension  $n$ de $M$ est paire, alors $H_n$ ne
  d\'epend plus de la m\'etrique.
Elle est, \`a un facteur constant pr\`es,
  la caract\'eristique d'Euler-Poincar\'e de $M$.\\
  M. Berger a d\'emontr\'e, voir  \cite{Berger}, que le  gradient 
de $H_4$,
  comme le gradient de $S$, d\'epend seulement du tenseur de
 courbure de Riemann  $R$
   et non pas de ses d\'eriv\'ees covariantes. Il a alors pos\'e la question si ce ph\'enom\`ene
reste vrai pour toutes les autres $H_{2k}$.
\par\noindent
  Le th\'eor\`eme suivant donne une r\'eponse affirmative \`a cette question:
\par\medskip
\noindent
  \begin{theorem}[\cite{Labbi9}]
   Soit $(M,g)$ une vari\'et\'e  riemannienne
compacte de dimension $n$. Pour tout entier  $k$, tel
 que $2\leq 2k\leq n$,
la fonctionnelle  $H_{2k}$ est  diff\'erentiable, et en  $g$ on a
  $$H_{2k}'h={1\over 2}<T_{2k},h>.$$
O\`u $T_{2k}$ est le tenseur d'Einstein-Lovelock d'ordre $2k$ d\'efini ci-dessus.
\end{theorem}
\par\noindent\medskip
Remarquons que  pour $k=1$, on a
$H_{2}'h={1\over 2}<T_2,h>={1\over 2}<{{\rm scal}\over 2} g- {\rm Ric},h>.$
Ceci n'est autre que la formule classique pour la courbure scalaire totale.
Aussi, pour $2k=n$, on obtient
$H_{n}'h={1\over 2}<T_n,h>=0.$
Ceci \'etant pr\'evisible car  $H_n$ ne d\'epend plus de la m\'etrique, comme l'affirme le th\'eor\`eme
de Gauss-Bonnet.\\
Notons que ce r\'esultat a \'et\'e premi\`erement \'etabli par Lovelock \cite{Lovelock} en
utilisant le calcul tensoriel classique. Cet article  \cite{Lovelock}
reste tr\`es peu
connu dans les milieux
math\'ematiques.
\par\medskip
\noindent
{\sc Esquisse de la d\'emonstration.}
On note tout d'abord que la d\'eriv\'ee directionnelle du tenseur de courbure de Riemann vu comme une forme
double sym\'etrique s'\'ecrit:
\begin{equation}\label{derivative:R}
R'h={-1\over 4}(D\tilde D+\tilde D D)(h)+{1\over 4}F_h(R).
\end{equation}
o\`u  $(D\tilde D+\tilde D D)$ est l'op\'erateur Hessien g\'en\'eralis\'e, voir chapitre 1,
et pour tout  $h\in {\cal C}^1$,  $F_h:{\cal C}\rightarrow
{\cal C}$ est un op\'erateur auto-adjoint qui agit par d\'erivations sur ${\cal C}$. En particulier,
$$
F_h(R)= h(R(x,y)z,u)
-h(R(x,y)u,z)+
h(R(z,u)x,y)-h(R(u,z)x,y).$$
Ensuite, on calcule la d\'eriv\'ee directionnelle de $h_{2k}$ en $g$.
On obtient,
$$h_{2k}'h={-1\over 2}<{c^{2k-1}\over (2k-1)!}R^k,h>+
  (-1)^n{k\over 4}(\delta\tilde\delta +\tilde\delta\delta )\biggl( *
  ( {g^{n-2k}\over (n-2k)!}R^{k-1}h)\biggr).$$
  O\`u $ (\delta\tilde\delta +\tilde\delta\delta )$ esl l'adjoint formel du
  Hessien $(D\tilde D+\tilde D D)$, voir chapitre 1.\\
Enfin, la preuve s'ach\`eve en utilisant le th\'eor\`eme de Stockes comme suit
\begin{equation*}
\begin{split}
H_{2k}'.h&=\int_M\biggl(h_{2k}'.h+{h_{2k}\over 2}{\rm
tr}_gh\biggr)\mu_g\\
&=-{1\over 2}<{c^{2k-1}\over (2k-1)!}R^k,h>+{h_{2k}\over 2}<g,h>\\
&={1\over 2}<h_{2k}g-{c^{2k-1}\over (2k-1)!}R^k,h>
={1\over 2}<T_{2k},h>.
\end{split}
\end{equation*}
\ppp
\subsection{Le probl\`eme de Yamabe g\'en\'eralis\'e}
Comme cons\'equence directe du th\'eor\`eme pr\'ec\'edent, on a d\'emontr\'e que:
\begin{proposition}[\cite{Labbi9}] Soit $(M,g)$ une vari\'et\'e  riemannienne compacte  de
  dimension
$n>2k$. La courbure de Gauss-Bonnet-Weyl    $h_{2k}$ est constante si et seulement si
 la m\'etrique $g$
est un point critique de la fonctionnelle  $H_{2k}$ restreinte \`a
l'ensemble ${\rm Conf}_0(g)$ des m\'etriques ponctuellement conformes \`a $g$
et ayant le m\^eme volume total.\end{proposition}
Il est alors naturel de se poser la question si pour  tout $k$ on a: \\
\par\medskip {\sl Dans  toute classe conforme d'une m\'etrique
riemannienne sur une vari\'et\'e  compacte, il existe une m\'etrique riemannienne
\`a courbure de
Gauss-Bonnet-Weyl   $h_{2k}$ constante.}\\
\par\medskip\noindent
Le probl\`eme suivant, dit $\sigma_k$-probl\`eme de Yamabe, est \'etroitement li\'e  \`a
la question pr\'ec\'edente. Elle est
 actuellement l'objet de recherches intensives, voir le compte rendu de 
Viacklovsky \cite{Viaclovsky}, et on peut l'\'enoncer comme suit:\\
\par\smallskip
{\sl Notons par $\sigma_k$ les fonctions sym\'etriques des valeurs propres du tenseur
de Schouten.
 Pour tout k, il existe dans chaque classe conforme d'une m\'etrique donn\'ee sur 
une vari\'et\'e compacte,
une m\'etrique
riemannienne \`a courbure $\sigma_k$ constante.}\\
\par\smallskip
Ces deux probl\`emes co\"{\i}ncident dans le cas conform\'ement plat si $k$ est pair, puisque
 les deux  courbures $h_{2k}$ et $\sigma_{2k}$ co\"{\i}ncident \`a un facteur constant
  pr\`es. Aussi, dans  le cas $k=1$, les deux probl\`emes co\"{\i}ncident avec
le c\'el\`ebre probl\`eme de
Yamabe.\\
Le $\sigma_k$-probl\`eme de Yamabe a \'et\'e r\'ecemment r\'esolu   pour $k>n/2$  par
Gursky et Viaclovsky
\cite{Gursky} en supposant que la m\'etrique est ``admissible''. Ensuite  Sheng, Trudinger et
Wang \cite{Trudinger} ont compl\'et\'e les cas  o\`u $2\leq k\leq n/2$ en imposant
en plus que l'\'equation en question soit variationnelle.
\\
Notons toutefois qu'auparavant ce m\^eme probl\`eme a \'et\'e r\'esolu dans le cas conform\'ement plat
par Li et Li \cite{Lili} et Guan et Wang \cite{Guan}.\\
Soulignons enfin que, contrairement \`a ce qu'on peut  comprendre de \cite{Trudinger},
les fonctions $h_{2k}$ sont en g\'en\'eral  diff\'erentes des fonctions sym\'etriques
des valeurs propres de l'op\'erateur de courbure. Ces derni\`eres sont au nombre de $n(n-1)/2$
qui est nettement plus grand que $n/2$. Aussi,  on peut se rendre compte de la diff\'erence
entre ces courbures en consid\'erant le cas  
d'une hypersurface de l'espace Euclidien ou m\^eme dans le cas conform\'ement plat.
 Toutefois, il serait utile d'\'etablir
des relations alg\'ebriques entre ces invariants.

\subsection{Vari\'et\'es d'Einstein g\'en\'eralis\'ees}
Les m\'etriques d'Einstein usuelles sont les points critiques de la fonctionnelle
de la courbure scalaire totale, une fois restreinte aux m\'etriques de volume 1. Par analogie,
en consid\'erant les points critiques de la fonctionnelle de la courbure de Gauss-Bonnet-Weyl
 totale
on obtient des vari\'et\'es d'Einstein g\'en\'eralis\'ees. Plus pr\'ecis\'ement:\\ 
Pour $2\leq 2k\leq n$, on dira que  $(M,g)$ est une vari\'et\'e  $(2k)$-Einstein si le tenseur
 d'Einstein-Lovelock d'ordre  $2k$ est
proportionel \`a la m\'etrique, i.e.
$$T_{2k}=\lambda g.$$
 Du fait que  les  tenseurs $T_{2k}$ sont \`a divergence nulle, la fonction $\lambda$
 est alors constante.
 Remarquons aussi que les vari\'et\'es 2-Einstein sont  les vari\'et\'es d'Einstein usuelles.
De plus, pour $n=2k$, on a $T_{2k}=0$, alors toute m\'etrique riemannienne sur une vari\'et\'e
de dimension $n$ est $n$-Einstein.
\par\medskip\noindent
 La classe des vari\'et\'es d'Einstein g\'en\'eralis\'ees d'ordre $ 2k$
 contient les vari\'et\'es \`a courbure sectionnelle constante et toutes les
 vari\'et\'es homog\`enes \`a isotropie irr\'eductible munies de leurs
m\'etriques riemanniennes naturelles.\\
Dans les lignes qui suivent, je vais pr\'esenter  quelques remarques sur ces m\'etriques.:\\
\begin{itemize}
\item  Soit $ M$ une vari\'et\'e riemannienne de dimension
$n$, le produit riemannien standard $M\times {\bf R}^q$ est tel que
 $T_{2k}=0$
pour $2k\geq n$, mais  $T_2$ n'est pas en g\'en\'eral identiquement nulle.
Cet exemple nous am\`ene \`a penser naturellement que la condition pour une m\'etrique
 d'\^etre d'Einstein
au sens usuel est plus forte. Cependant, ceci n'est pas toujours le cas comme le montre
le contre exemple suivant:
\item Soit  $M$ une vari\'et\'e riemannienne de dimension 4  Ricci-plate mais pas
plate, par exemple une surface  K3 munie  de la m\'etrique de Calabi-Yau. Si
$T^q$ est le tore plat, alors le produit riemannien $M\times T^q$ est d'Einstein
au sens usuel. En revanche celui-ci n'est pas 4-Einstein.
\item Rappelons que, voir chapitre 1, le tenseur de Gauss-Kronecker admet la d\'ecomposition
suivante
$$R^q=\omega_{2q}+g\omega_{2q-1}+...+g^{2q-1}\omega_{1}+g^{2q}\omega_0$$
On montre  \cite{Labbi7} qu'une m\'etrique est $(2q)$-Einstein si et seulement si
la composante
$\omega_1$ est identiquement nulle. Ceci g\'en\'eralise un r\'esultat similaire
pour les m\'etriques d'Einstein usuelles.\\
\end{itemize}

\par\medskip\noindent
\section{Les $(p,q)$-courbures}

\subsection{Introduction}
La  m\'etrique $g$ et le tenseur de courbure de  Riemann  $R$ satisfont la premi\`ere et
la deuxi\`eme identit\'e
de  Bianchi, alors il en est de m\^eme pour tous les produits  $g^pR^q$. Par dualit\'e,
les tenseurs  $*(g^pR^q)$ satisfont eux aussi  la premi\`ere identit\'e
de  Bianchi et sont tous \`a divergence nulle.\\
Ces tenseurs de $(p,q)$-courbure sont compl\`etement d\'etermin\'es \`a partir de leurs courbures
sectionnelles. Ils englobent plusieurs
courbures bien connues, y compris,   les courbures: scalaire,
sectionnelle, de Gauss-Kronecker,   toutes les $p$-courbures,
tous
les invariants
 de courbure de la formule des tubes de Weyl, ainsi que
 les courbures d'Einstein-Lovelock ...
\\
Commen\c cons tout d'abord par la d\'efinition:
\begin{definition} La $(p,q)$-courbure, not\'ee  $s_{(p,q)}$,
  pour $1\leq q\leq {n\over 2}$ et $0\leq p\leq n-2q$, est la courbure sectionnelle
du  $(p,q)$-tenseur de courbure suivant:
  \begin{equation}\label{spq:def}
  R_{(p,q)}={1\over (n-2q-p)!}*\bigl( g^{n-2q-p}R^q\bigr).
  \end{equation}
 Pr\'ecis\'ement, pour un $p$-plan tangent, $s_{(p,q)}(P)$ est la courbure sectionnelle
du tenseur
  ${1\over (n-2q-p)!}g^{n-2q-p}R^q$ dans la direction du $(n-p)$-plan suppl\'ementaire et
   orthogonal \`a $P$.
\end{definition}
Remarquons que les tenseurs  $R_{(p,q)}$ satisfont la premi\`ere identit\'e de Bianchi
 mais en g\'en\'eral pas la deuxi\`eme. En revanche,  ils sont toujours \`a divergence nulle. En effet,
$$\delta R_{(p,q)}=\delta {1\over (n-2q-p)!}*\bigl( g^{n-2q-p}R^q\bigr)
=*D{1\over (n-2q-p)!}\bigl( g^{n-2q-p}R^q\bigr)=0.$$
En particulier, on a montr\'e que:\\
\begin{proposition}[Th\'eor\`eme de Schur \cite{Labbi9}]
 Soit $p,q\geq 1$.
 Si en tout point  $m\in M$, la  $(p,q)$-courbure est constante,
 alors elle est constante.\end{proposition}
Ces courbures g\'en\'eralisent plusieurs notions de courbures bien connues.
Notons que pour  $q=1$, on a  $s_{(p,1)}=s_p$, o\`u $s_p$ est la $p$-courbure
que j'avais \'etudi\'e dans ma th\`ese, voir ci-dessous.
 En particulier,  $s_{(0,1)}$ est la moiti\'e de la courbure scalaire,
 et   $s_{(n-2,1)}$ co\"{\i}ncide avec la courbure sectionnelle de  $(M,g)$.\\
Aussi, pour  $p=0$ et $2q=n$, $s_{(0,{n\over 2})}=*R^{n/2}$ est,
\`a un facteur constant
pr\`es, la courbure de  Killing-Lipschitz. Plus g\'en\'eralement,
la courbure $s_{(n-2q,q)}(P)$ est, \`a un facteur constant pr\`es, la courbure de
 Killing-Lipschitz
de $P^{\bot}$. Cette derni\`ere  n'est autre que  la  $(2p)$-courbure sectionnelle d\'efinie par
Thorpe, voir \cite{Thorpe}.\\
Pour $p=0$, $s_{(0,q)}=*{1\over (n-2q)!}g^{n-2q}R^q={1\over (2q)!}c^{2q}R^q$
 est la courbure de Gauss-Bonnet-Weyl introduite dans le chapitre 2.
 \\
 Finalement, pour  $p=1$, $s_{(1,q)}$ est la courbure sectionnelle du tenseur
 d'Einstein-Lovelock \'etudi\'ee aussi  dans le chapitre 2.\\
\par\medskip\noindent

\subsection{La $p$-courbure}
Dans cette section, on va s'int\'eresser au cas particulier $q=1$. On retrouve alors
la $p$-courbure.\\
Rappelons que la $p$-courbure, not\'ee $s_p$, est d\'efinie pour $0\leq p\leq n-2$
en \'etant la courbure sectionnelle du tenseur
 $${1\over (n-2-p)!}*(g^{n-2-p}R).$$
Pour un $p$-plan tangent en  $m\in M$, $s_p(P)$ est la moiti\'e de la courbure scalaire
en $m$
de la sous vari\'et\'e totalement g\'eod\'esique  ${\rm exp}_m\, P^{\bot}$ de $M$. Pour $p=0$,
on retrouve la courbure scalaire usuelle, et pour  $p=n-2$ elle co\"{\i}ncide avec la
courbure sectionnelle.
\par\medskip\noindent
 En utilisant les r\'esultats du chapitre 1, on a d\'emontr\'e ais\'ement
le th\'eor\`eme ci-dessous. Celui-ci  \'etabli une caract\'erisation g\'eom\'etrique
des m\'etriques \`a courbure
sectionnelle constante, des m\'etriques conform\'ement plates \`a courbure scalaire
constante et aussi des m\'etriques d'Einstein. Signalons que des r\'esultats semblables ont
\'et\'e d'abord d\'emontr\'es
 dans  ma th\`ese et ensuite par Chen-Dillen-Verstraelen-Vrancken dans \cite{Shin}
 en utilisant des calculs assez longs!

\begin{theorem}[\cite{Labbi7}]\label{sp+sp:sp-sp}\begin{enumerate}
\item  Pour tout   $2\leq p\leq n-2$, la $p$-courbure est
constante si et seulement si  $(M,g)$ est \`a courbure sectionnelle constante.
\item Pour tout $1\leq p\leq n-1$,   $(M,g)$ est une vari\'et\'e
 d'Einstein si et seulement si la fonction $P\rightarrow s_p(P)-s_{n-p}(P^{\bot})
 =\lambda$ est
 constante. De plus, si tel est le cas on a  $\lambda={n-2p\over 2n}c^2R.$
\item Pour tout   $2\leq p\leq n-2$ et $p\not = {n\over 2}$, la fonction
  $P\rightarrow s_p(P)+s_{n-p}(P^{\bot})=\lambda$ est constante si et seulement si
$(M,g)$
est \`a courbure sectionnelle constante.
De plus, si tel est le cas on a  $\lambda={2p(p-1)+(n-2p)(n-1)\over 2n(n-1)}c^2R.$
\item Soit  $n=2p$, alors  $(M,g)$ est conform\'ement plate \`a courbure scalaire constante
si et seulement si  la fonction
  $P\rightarrow s_p(P)+s_{p}(P^{\bot})=\lambda$ est constante.
De plus, dans ce cas on aura  $\lambda={n-2\over4(n-1)} c^2R.$
\end{enumerate}
\end{theorem}
{\sc Id\'ee de d\'emonstration}. Les conditions pr\'ec\'edentes se lisent au niveau
des tenseurs de courbure correspondants comme  \'etant proportionnels \`a une
certaine puissance de la m\'etrique. On ach\`eve  la d\'emonstration,  tout simplement,
apr\`es une
division des deux membres de l'\'equation obtenue (au niveau des
tenseurs)
par une puissance ad\'equate  de la m\'etrique $g$. Et ceci gr\^ace \`a la proposition 1
 du chapitre 1.

\subsection{Propri\'et\'es  des $(p,q)$-courbures}
Dans cette section on pr\'esente les g\'en\'eralisations des r\'esultats pr\'ec\'edents
au cas  des
$(p,q)$-courbures. Les d\'emonstrations sont identiques au cas de la $p$-courbure.\\
\medskip
 Le r\'esultat suivant caract\'erise les m\'etriques \`a $(p,q)$-courbure
 constante.
\begin{proposition}[\cite{Labbi7}]
\begin{enumerate}
\item Pour tout $(p,q)$ tel que  $2q\leq p\leq n-2q$,
 la $(p,q)$-courbure
 $s_{(p,q)}\equiv \lambda$ est constante si et seulement si la courbure sectionnelle
du tenseur de Gauss-Kronecker
 $R^q$ est constante et \'egale \`a  ${\lambda(2q)!(n-p-2q)!\over (n-p)!}$.
 \item Pour tout  $(p,q)$ tel que  $p<2q$, la  $(p,q)$-courbure
 $s_{(p,q)}\equiv c$ est constante si et seulement si le tenseur
  $c^{2q-p}(R^q)$ est proportionel \`a la m\'etrique, i.e.
   $c^{2q-p}(R^q)=\lambda g^p$.
   \end{enumerate}
   \end{proposition}

Rappelons qu'une m\'etrique est dite d'Einstein si son tenseur de Ricci $cR$  lui est
 proportionnel. La condition 2 pr\'ec\'edente peut \^etre consid\'er\'ee alors 
comme une g\'en\'eralisation de la condition d'Einstein.\\
Le r\'esultat suivant sp\'ecifie cette condition et g\'en\'eralise un r\'esultat
similaire pour les m\'etriques d'Einstein:
\begin{proposition}[\cite{Labbi7}]\label{prop1}
 Pour $1\leq p<2q$,  le tenseur  $c^{2q-p}(R^q)$ est proportionel \`a
la m\'etrique $g^p$ si et seulement si
$$\omega_i=0 \, {\text { pour}}\,  1\leq i\leq {\min\{p,n-p\}},$$
o\`u les $\omega_i$ sont les composantes de $R^q$ suivant la d\'ecomposition orthogonale
\ref{ort:decom} donn\'ees par 
$R^q=\sum_{i=0}^{2q}g^{2q-i}\omega_i$.\end{proposition}
\medskip
Finalement, on a caract\'eris\'e g\'eom\'etriquement 
 l'annulation des diff\'erentes composantes $\omega_i$ dans la d\'ecomposition orthogonale
de $R^q$ et ceci en utilisant les $(p,q)$-courbures.
\begin{proposition}[\cite{Labbi7}]
Soit
 $2q\leq r\leq n-2q$, $n\not =2r$ et
$R^q=\sum_{i=0}^{2q}g^{2q-i}\omega_i$,
alors
\begin{enumerate}
\item La fonction $P\rightarrow s_{(r,q)}(P)- s_{(n-r,q)}(P^\bot) \equiv
 \lambda$ est constante
si et seulement si $\omega_i=0$ pour $1\leq i\leq 2q-1$ et
 $\bigl({(n-r)!\over (n-2q-r)!}-{r!\over (r-2q)!}\bigr)\omega_0=\lambda$.
\item La fonction $P\rightarrow s_{(r,q)}(P)+ s_{(n-r,q)}(P^\bot) \equiv
 \lambda$ est constante
si et seulement si $\omega_i=0$  pour $1\leq i\leq 2q$ et $\bigl({(n-r)!\over
 (n-2q-r)!}+{r!\over (r-2q)!}\bigr)\omega_0=\lambda$. Autrement
 dit,
  $R^q$ est \`a courbure sectionnelle constante.
  \end{enumerate}
\end{proposition}
Les cas restants sont discut\'es dans la proposition suivante:
\begin{proposition}[\cite{Labbi7}]
Soit  $2q\leq r\leq n-2q$ et $n =2r$, alors
\begin{enumerate}
\item La fonction $P\rightarrow s_{(r,q)}(P)- s_{(r,q)}(P^\bot) \equiv
\lambda$ est constante
si et seulement si $\omega_i=0$ pour les $i$ impairs tels que
 $1\leq i\leq 2q-1$
et $\lambda=0$.
\item La fonction $P\rightarrow s_{(r,q)}(P)+ s_{(r,q)}(P^\bot) \equiv
 \lambda$ est constante
si et seulement si $\omega_i=0$ pour les $i$ pairs tels que
 $2\leq i\leq 2q$ et
 $2{r!\over (r-2q)!}\omega_0=\lambda$.
\end{enumerate}
\end{proposition}

\subsection{Formule d'Avez g\'en\'eralis\'ee}
Rappelons que la formule d'Avez \cite{Avez} montre que l'int\'egrand de Gauss-Bonnet
en dimension 4 est
donn\'e par $h_4=\sum_{r=0}^{2}{(-1)^{r}\over (r!)^2}|c^rR|^2$. Cette formule a servi en particulier
\`a d\'emontrer que les vari\'et\'es d'Einstein compactes et de dimension 4 ont leurs caract\'eristiques
d'Euler-Poincar\'e non-n\'egative.\\
On a g\'en\'eralis\'e le r\'esultat pr\'ec\'edent aux dimensions sup\'erieures comme suit:\\
\begin{theorem}[\cite{Labbi7}] Soit $n=2p$ et $\omega, \theta\in {\cal C}_1^p$, alors
$$*(\omega\theta)=\sum_{r=0}^p{(-1)^{r+p}\over (r!)^2}<c^r\omega,c^r\theta>.$$
En particulier, si $n=4q$, alors l'int\'egrand de  Gauss-Bonnet est donn\'e par
$$h_{4q}=\sum_{r=0}^{2q}{(-1)^{r}\over (r!)^2}|c^rR^q|^2.$$
\end{theorem}
\par\noindent\medskip
Le r\'esultat suivant est du m\^eme type:
\begin{theorem}[\cite{Labbi7}]\label{thm1} Suivant la d\'ecomposition orthogonale  \ref{ort:decom},
soit
$\omega=\sum_{i=0}^{n-p}g^{n-p-i}\omega_i\in {\cal C}_1^{n-p}$ et
$\theta=\sum_{i=0}^p g^{p-i}\theta_i\in {\cal C}_1^p$ , alors
$$*(\omega\theta)=\sum_{r=0}^{\min\{p,n-p\}}(-1)^r(n-2r)!<\omega_i,\theta_i>.$$
En particulier, pour tout $q$ tel que $n\geq 4q$ on a
$$h_{4q}={1\over (n-4q)!}\sum^{2q}_{i=0}(-1)^i
(n-2i)!<(R^q)_i,(R^q)_i>.$$
\end{theorem}
Le cas $q=1$ est particuli\`erement int\'eressant.  En effet, si $R=\omega_2+g\omega_1
+g^2\omega_0$ est la d\'ecomposition standard du tenseur de courbure de Riemann. Le
th\'eor\`eme  pr\'ec\'edent affirme qu'en dimension $\geq 4$ on a
$$(n-4)!h_4=n!|\omega_0|^2-(n-2)!|\omega_1|^2+(n-4)!|\omega_2|^2.$$
Par cons\'equent, on  a pu d\'emontrer que:
\begin{theorem}[\cite{Labbi7, Labbi8}]\begin{enumerate}
\item Pour une vari\'et\'e d'Einstein  $(M,g)$ de dimension  $n\geq 4$ on a
$h_4\geq 0$. De plus  $h_4\equiv 0$ si et seulement si $(M,g)$ est plate.
\item Si une vari\'et\'e conform\'ement plate $(M,g)$ est  \`a courbure scalaire
nulle et de dimension
  $n\geq 4$,
alors $h_4\leq 0$. De plus,  $h_4\equiv 0$ si et seulement si
 $(M,g)$ est plate.
\end{enumerate}
\end{theorem}
On en d\'eduit alors  une
obstruction g\'eom\'etrique \`a l'existence des m\'etriques d'Einstein
ainsi que pour les m\'etriques conform\'ement plates \`a courbure scalaire nulle.\\
Rappelons qu'en dimension $\geq 5$ on ne  connait pas d'obstructions topologiques
\`a l'existence
de m\'etriques d'Einstein. Jusqu'ici, les seules obstructions connues  sont celles
obtenues en imposant la condition suppl\'ementaire de  positivit\'e de la 
 constante d'Einstein. Celle-ci force la courbure scalaire et
la courbure de Ricci  d'\^etres
positives,  et on sait qu'on a des r\'estrictions topologiques dans ces cas.\\
L'obstruction g\'eom\'etrique pr\'ec\'edente \`a l'avantage 
d'\^etre ind\'ependante
du signe de la constante d'Einstein. C'est d'ailleurs ce r\'esultat qui a
motiv\'e notre \'etude de la positivit\'e de la seconde courbure de
Gauss-Bonnet-Weyl.
\par\noindent\medskip
Le th\'eor\`eme  pr\'ec\'edent se g\'en\'eralise naturellement
 pour les courbures sup\'erieures comme suit:
\begin{theorem} Soit $(M,g)$ une vari\'et\'e riemannienne de dimension
  $n\geq 4q$ et telle que le tenseur $c(R^q)$ est proportionnel \`a
 la m\'etrique $g^{2q-1}$
alors $h_{4q}\geq 0$. De plus dans ce cas
$h_{4q}\equiv 0$  si et seulement si  $(M,g)$ est $q$-plate.\\
En particulier, une vari\'et\'e compacte de dimension $n=4q$ telle que son tenseur $c(R^q)$ 
est proportionnel \`a la m\'etrique $g^{2q-1}$ est \`a caract\'er\'estique d'Euler-Poincar\'e
non-n\'egative. De plus, Elle est nulle si et seulement si la vari\'et\'e est $q$-plate.
\end{theorem}
Rappelons que $q$-plate veut dire que la courbure sectionnelle de $R^q$ est identiquement nulle.
Notons que ce r\'esultat n'est pas publi\'e dans mes travaux. Il se d\'emontre directement
\`a partir du th\'eor\`eme \ref{thm1} et de la proposition \ref{prop1}.\\
Notons \'egalement que cette condition est plus forte que la condition $(2q)$-Einstein d\'efinie 
dans le chapitre 2.\\
\section{Les courbures de Weitzenb\"{o}ck}
L'op\'erateur de courbure de Weitzenb\"{o}ck, not\'e ici par  ${\cal N}$,
est le terme d'ordre z\'ero (i.e.
d\'ependant lin\'eairement de la courbure)
 dans la 
c\'el\`ebre formule de Weitzenb\"{o}ck.
Cette formule exprime le Laplacien $\Delta$ sur les formes diff\'erentielles
 en termes de la connexion
de Levi civita $\nabla$, pr\'ecis\'ement:\\
$$\Delta=\nabla^*\nabla+{\cal N}.$$
Cette formule est importante dans l'\'etude des int\'eractions entre la 
g\'eom\'etrie  et la topologie d'une vari\'et\'e. En effet, il y a une
m\'ethode qui remonte \`a Bochner et connue sous le nom de
th\'eor\`emes d'annulations, consistant \`a 
d\'emontrer
l'annulation des nombres de Betti d'une vari\'et\'e riemannienne ayant
 certaines positivit\'es
de courbure entra\^{\i}nant la positivit\'e de ${\cal N}$.
 Cette m\'ethode s'applique  essentiellement sur les vari\'et\'es compactes.\\
L'op\'erateur ${\cal N}$ envoie les $p$-formes sur elles m\^emes et il
 est auto-adjoint.
Ceci permet, par dualit\'e,
de le voir comme une forme double.\\
Le probl\'eme avec ${\cal N}$ est qu'il s'est toujours impos\'e comme
  une expression compliqu\'ee de la courbure
et donc difficilement manipulable. Diff\'erentes simplifications de cet  op\'erateur existent.
 Je cite par exemple
le formalisme de Clifford dans \cite{Lawson-Michelsohn, Labbi1}, 
le travail de Gallot-Meyer dans 
\cite{Gallot-Meyer},.... \\
Dans notre contribution \cite{Labbi11}, on montre
 une formule  simple pour $N$ qu'on utilise par la suite pour
 \'etudier des
propri\'et\'es g\'eom\'etriques de cette courbure. Signalons que
 cette formule  de
${\cal N}$ a \'et\'e auparavant \'etablie par Bourguignon dans 
\cite{Bourguignon},  proposition 8.6. \`A ma connaissance  cette formule  
 reste malheureusement m\'econnue et
pas utilis\'ee. Pour ma part, je n'ai pris connaissance de celle-ci
qu'une fois ma preuve achev\'ee. Notre d\'emonstration est compl\`etement diff\'erente
et est directe. \\
\subsection{Courbures sectionnelles de Weitzenb\"{o}ck}
Un calcul direct sans difficult\'es montre que la courbure
sectionnelle du tenseur ${\cal N}$
est donn\'ee par
\begin{equation*}
{\cal N}(e_{i_1}\wedge e_{i_2}\wedge ...\wedge e_{i_p};
e_{i_1}\wedge e_{i_2}\wedge ...\wedge e_{i_p})
=\sum_{j=1}^p
\sum_{k=p+1}^n K(e_{i_j},e_{i_k}). 
\end{equation*}
O\`u $(e_{i_1},...,e_{i_n})$ est une base orthonorm\'ee de vecteurs
 tangents et $K$ d\'esigne la courbure sectionnelle du tenseur de courbure de Riemann $R$.\\
D'autre part, on note que
\begin{equation*}
\begin{split}
\sum_{i=1}^p\sum_{j=p+1}^n R(e_i\wedge e_j &,e_i\wedge e_j)=
\sum_{i=1}^p
cR(e_i,e_i)
-\sum_{i,j=1}^pR(e_i\wedge e_j,e_i\wedge e_j)\\
&=\{\frac{g^{p-1}}{(p-1)!}cR -2\frac{g^{p-2}}{(p-2)!}
R\}(e_{1}\wedge...\wedge e_{p},
e_{1}\wedge...\wedge e_{p}).
\end{split}
\end{equation*}
Il en  r\'esulte alors imm\'ediatement  que la courbure sectionnelle de
 Weitzenb\"{o}ck co\"incide avec la courbure sectionnelle de la forme double $\{ \frac{gcR}{(p-1)}-2R\}
\frac{g^{p-2}}{(p-2)!}$.
 
\subsection{Structures de courbures de Weitzenb\"{o}ck}
Notons que la forme double $\{\frac{gcR}{(p-1)}-2R\}\frac{g^{p-2}}{(p-2)!}$ satisfait la premi\`ere
identit\'e de Bianchi,  elle a en outre la m\^eme courbure sectionnelle que la structure de courbure de
Weitzenb\"{o}ck.
Il suffit alors de d\'emontrer que cette structure satisfait elle aussi la premi\`ere
identit\'e de Bianchi (voir \cite{Labbi11})
 pour avoir une d\'emonstration compl\`ete de la formule suivante:
\begin{proposition}[\cite{Bourguignon,Labbi11}]
Pour tout $p$ tel que $2\leq  p\leq n-2$, le tenseur de courbure de la
 formule de Weitzenb\"{o}ck
n'est autre que
\begin{equation}
{\cal N}=\{\frac{gcR}{(p-1)}-2R\}\frac{g^{p-2}}{(p-2)!}.
\end{equation}
\end{proposition}
Rappelons que pour $p=1$, la forme double ${\cal N}$ n'est autre que le tenseur de Ricci.
\subsection{Propri\'et\'es  g\'eom\'etriques}
Pour pr\'eciser le degr\'e $p$ de la forme double ${\cal N}$, on \'ecrira  ${\cal N}_p$.\\
Le r\'esultat suivant peut \^etre d\'emontr\'e ais\'ement, soit  en utilisant la formule explicite
de l'op\'erateur de Hodge \'etablie dans le premier chapitre soit  tout simplement en remarquant
que les deux formes doubles satisfont la premi\`ere identit\'e de Bianchi et ont la m\^eme
courbure sectionnelle:
\begin{proposition}[\cite{Labbi11}]
Pour tout $p$ tel que  $ 2\leq p\leq n-2$, on a
$$* {\cal N}_p={\cal N}_{n-p}.$$
En particulier, si $n=2p$ on obtient $*{\cal N}_p={\cal N}_{p}.$
\end{proposition}
Il est clair que ${\cal N}_p$ est divisible par $g^{p-2}$ et donc sa d\'ecomposition
orthogonale sera r\'eduite comme suit:
\begin{proposition}[\cite{Labbi11}]
Soit $R=\omega_2+g\omega_1+g^2\omega_0$ la d\'ecomposition standard du tenseur de courbure de Riemann,
alors le tenseur ${\cal N}_p$ se d\'ecompose suivant \ref{ort:decom} comme suit:
$${\cal N}_p=-g^{p-2}\omega_2+g^{p-1}\{\frac{n-2p}{2(p-1)}\omega_1\}+
g^p\{\frac{n-p}{p-1}\omega_0\}.$$
Ceci \'etant bien s\^ur pour tout $2\leq p\leq n-2$.
\end{proposition}
Rappelons qu'une forme double $(p,p)$ satisfaisant la premi\`ere identit\'e de Bianchi
est \`a  courbure sectionnelle constante si et seulement si elle est proportionelle \`a $g^p$.
Le corollaire suivant est  alors imm\'ediat  en utilisant soit la proposition pr\'ec\'edente
soit  la proposition 1 du chapitre 1.: \\
\begin{corollary}[\cite{Labbi11}] Pour $2\leq p\leq n-2$, une vari\'et\'e riemannienne $(M,g)$ de dimension
 $n$ est \`a  courbure
sectionnelle de Weitzenb\"{o}ck d'ordre $p$ constante si et seulement si elle est soit 
  \`a courbure sectionnelle constante soit
 conform\'ement plate de dimension $n=2p$.
\end{corollary}
Signalons qu'un r\'esultat du m\^eme type a \'et\'e d\'emontr\'e dans \cite{Tachibana,Nasu1}.
Les r\'esultats pr\'ec\'edents sont alg\'ebriques, i.e. valables en tout point
de la vari\'et\'e. En particulier, on en d\'eduit un th\'eor\`eme de  Schur pour ces courbures.\\
M\^eme si l'\'etude de la positivit\'e des courbures sera discut\'ee dans la deuxi\`eme partie
de ce m\'emoire, je me permettrai, exceptionellement ici, de mentionner quelques remarques
sur la positivit\'e de ${\cal N}_p$.\\
La positivit\'e de la forme double ${\cal N}_p$ est bien \'evidement 
tr\`es importante dans l'\'etude des
int\'eractions entre topologie et g\'eom\'etrie. Il serait donc tr\`es 
int\'eressant de comprendre
cette condition et en particulier, d'\'etablir des liens avec la positivit\'e des autres courbures.\\

Faut-il le rappeler ici que la positivit\'e de ${\cal N}_p$ est strictement plus forte que celle 
da sa courbure sectionnelle (penser \`a la positivit\'e de l'op\'erateur
de courbure de Riemann et de sa courbure sectionnelle).\\
On sait d\'ej\`a que la positivit\'e de l'op\'erateur de courbure de Riemann entra\^{\i}ne celle de 
${\cal N}_p$ pour tout $p$, c'est le th\'eor\`eme de Gallot-Meyer \cite{Gallot-Meyer}. Voir 
aussi \cite{Lawson-Michelsohn} pour une preuve simplifi\'ee. Notons aussi qu'une autre d\'emonstration
simplifi\'ee est possible  \`a partir de la formule de  ${\cal N}_p$ ci-dessus.\\
Dans \cite{Labbi11}  on \'etudie les liens avec la positivit\'e de la $p$-courbure et celle
de la courbure isotrope. En particulier, on  montre que:
\begin{proposition}[\cite{Labbi11}] Pour tout $2\leq p\leq n-2$, la contraction
compl\`ete de la forme double ${\cal N}_p$ est donn\'ee par
$$c^p({\cal N}_p)=\frac{p(n-2)!}{(n-p-1)!}c^2R$$
En particulier, la positivit\'e de ${\cal N}_p$ entra\^{\i}ne la positivit\'e de la courbure
scalaire.
\end{proposition}
D'une mani\`ere similaire, en prenant  des contractions jusqu'\`a l'ordre  $p-1$,
 on montre que la positivit\'e de
${\cal N}_p$ entra\^{\i}ne une condition de pincement sur les valeurs propres de la
 courbure de Ricci.\\
Dans le cas conform\'ement plat, on \'etablit la relation suivante entre le tenseur de
$p$-courbure
et ${\cal N}_p$:
\begin{proposition}[\cite{Labbi11}] Soit $(M,g)$ une vari\'et\'e conform\'ement
 plate de dimension $n$,
et $2\leq p\leq n-2$. Alors le tenseur ${\cal N}_p$ est donn\'e \`a partir du tenseur
de $p$-courbure comme suit:
$${\cal N}_{\frac{n+p}{2}}=\frac{p!}{n-p-1}g^{\frac{n-p}{2}}
\{*(\frac{1}{(n-p-2)!}g^{n-p-2}R)\}. $$
\end{proposition}
Rappelons que dans le conform\'ement plat, la positivit\'e de la $p$-courbure \'equivaut
\`a celle de son tenseur de $p$-courbure.
En particulier, la positivit\'e de la $p$-courbure entra\^{\i}ne alors la positivit\'e de
${\cal N}_{\frac{n+p}{2}}$. Par cons\'equent, on obtient une deuxi\`eme d\'emonstration 
plus rapide
de notre  r\'esultat sur l'annulation des nombres de Betti des vari\'et\'es conform\'ement
plates \`a $p$-courbure positive \cite{Labbi1}.
\newpage
\part{Variations sur diff\'erentes notions de positivit\'e de courbures}
\section{Positivit\'e de la $p$-courbure} 
L'\'etude de la positivit\'e de
la $p$-courbure a \'et\'e le sujet auquel \'etait consacr\'e  mon  premier travail  
de recherche. Cette courbure est une  extension de la courbure scalaire 
 propos\'ee par
M. Gromov et co\"{\i}ncide avec la courbure $p$-vectorielle de \'E. Cartan, 
voir \cite{Cartan,Nasu1}.
Elle est  une fonction d\'efinie
 sur la $p$-grassmannienne de la vari\'et\'e
en question, qui associe \`a tout $p$-plan $P$, la moyenne de la courbure sectionnelle dans la
direction du $(n-p)$-plan orthogonal. Pour $p=0$, ( resp. pour p=$n-2$), 
$n$ \'etant la dimension de la vari\'et\'e en question, on retrouve la courbure scalaire (resp.
la courbure sectionnelle). La positivit\'e de la $p$-courbure entra\^{\i}ne celle de la $(p-1)$-courbure,
pour tout $1\leq p\leq n-2$. En particulier, la positivit\'e de la $p$-courbure est 
interm\'ediaire
entre la positivit\'e des  courbures scalaire et sectionnelle. \\
En ce qui concerne 
 les probl\`emes portant
sur la courbure scalaire positive,  on a maintenant des r\'eponses tr\`es satisfaisantes
 suite aux
 travaux de Gromov-Lawson \cite{Gromov-Lawson}, Schoen-Yau, Stolz, ...
 Ceci n'\'etant absolument  pas  le cas pour la courbure sectionnelle positive. On ignore 
 toujours
si par exemple,  le produit $S^2\times S^2$ poss\`ede ou non une m\'etrique riemannienne \`a courbure
sectionnelle positive (conjecture de Hopf).\\
Afin de mieux situer mes travaux sur la positivit\'e, je commencerai tout d'abord  par un 
 rappel de mes 
r\'esultats de th\`ese. Signalons toutefois que les r\'esultats sur le groupe fondamental
 et ceux sur la courbure d'Einstein sont  post\'erieurs \`a ma th\`ese. J'ai
 d\'emontr\'e les r\'esultats suivants, qui g\'en\'eralisent des 
r\'esultats d\^us respectivement \`a Gromov-Lawson \cite{Gromov-Lawson}, Lawson-Yau 
\cite{Lawson-Yau}, Bourguignon \cite{Bourguignon} et Lafontaine \cite{Lafontaine-conf}.:\\
\begin{theorem}[\cite{Labbi3}] Si deux vari\'et\'es compactes de dimension $n\geq p+3$ 
admettent des m\'etriques
\`a $p$-courbure positive, il en est de m\^eme de leur somme connexe. Plus g\'en\'eralement,
si une vari\'et\'e compacte admet une m\'etrique \`a $p$-courbure positive, alors il en est de 
m\^eme
pour toute vari\'et\'e obtenue \`a partir de celle-ci par chirurgies de codimension
 $\geq p+3$.\end{theorem}
Comme cons\'equence de ce th\'eor\`eme et en utilisant des techniques de topologie alg\'ebrique,
 sp\'ecialement le cobordisme spinoriel de Gromov-Lawson, j'ai d\'emontr\'e des r\'esultats 
d'existence pour la 1-courbure positive  comme suit:
\begin{theorem}[\cite{Labbi3}]\begin{enumerate}
\item Une vari\'et\'e compacte  $2$-connexe de  dimension 
$\geq 7$ admet une m\'etrique riemannienne \`a $1$-courbure positive  si et seulement si 
elle admet une m\'etrique riemannienne \`a courbure scalaire positive. En particulier, toute 
vari\'et\'e compacte de dimension $7$ admet une m\'etrique \`a $1$-courbure
positive.
\item 
Toute vari\'et\'e compacte, simplement connexe, non-spinorielle  de dimension $\geq 7$
 et telle que son second groupe d'homotopie satisfait 
$\pi_2(V) \cong Z_2$ admet une m\'etrique riemannienne \`a $1$-courbure positive.
\end{enumerate}
\end{theorem}
Ensuite, j'ai utilis\'e \`a la mani\`ere de Lawson-Yau \cite{Lawson-Yau} les actions des groupes 
de Lie pour faire appara\^{\i}tre de la $p$-courbure positive. Le r\'esultat est le suivant:
\begin{theorem}[\cite{Labbi2}] Si une vari\'et\'e compacte admet une action effective d'un groupe de 
Lie compact simple de rang au moins $p+1$, alors cette vari\'et\'e porte une m\'etrique
\`a $p$-courbure positive.
\end{theorem}
En ce qui concerne les  obstructions, j'ai d\'emontr\'e un th\'eor\`eme d'annulation des nombres
de Betti pour les vari\'et\'es conform\'ement
plates  comme suit:
\begin{theorem}[\cite{Labbi1}]
\begin{enumerate}
\item Soit   $(M,g)$ une vari\'et\'e  compacte conform\'ement plate de dimension $n$  \`a
$p$-courbure positive alors les nombres de Betti sont nuls du degr\'e $\frac{n-p}{2}$
au degr\'e $\frac{n+p}{2}$.
\item Soit $(M,g)$ une vari\'et\'e  compacte conform\'ement plate de dimension $n$
\`a $p$-courbure   non-n\'egative  et telle que la cohomologie r\'eelle de degr\'e
$\frac{n\pm p}{2}$ est non nulle, alors   $(M,g)$ est soit  plate soit  isom\'etrique
\`a un quotient compact de 
$S^{\frac{n+p}{2}} \times H^{\frac{n-p}{2}}$.
\end{enumerate}
\end{theorem}
Pour d\'emontrer ces derniers r\'esultats, j'ai bien entendu utilis\'e la formule 
 de Weitzenb\"{o}ck mais reformul\'ee en termes  des alg\`ebres de Clifford selon les id\'ees
d\'evolop\'ees par Lawson-Michelsohn \cite{Lawson-Michelsohn}.\\
Par ailleurs, du fait que le produit d'une sph\`ere de dimension $\geq p+2$
 avec une vari\'et\'e compacte quelconque admet une m\'etrique \`a $p$-courbure positive, on
en d\'eduit que tout groupe de pr\'esentation finie peut \^etre r\'ealis\'e comme le groupe
fondamental d'une vari\'et\'e de dimension $\geq p+6$ et \`a $p$-courbure positive. En utilisant
notre r\'esultat pr\'ec\'edent sur les chirurgies, on a pu dire mieux comme suit:
\begin{theorem}[\cite{Labbi5}]
Soit $G$ un groupe \`a pr\'esentation finie et $p\geq 0$. Alors, pour tout $n\geq p+4$, il
existe une vari\'et\'e  compacte de dimension $n$ \`a $p$-courbure
 positive et telle que $\pi_1(M)=G$.
\end{theorem}

\section{ Courbure d'Einstein positive}
La courbure d'Einstein est une fonction  d\'efinie sur le fibr\'e tangent unitaire de
la vari\'et\'e en question.
Elle est obtenue \`a partir du c\'el\`ebre tenseur d'Einstein de la m\^eme  mani\`ere
que l'on
 d\'efinit la courbure de Ricci \`a partir du  tenseur de Ricci. Elle co\"{\i}ncide avec la 
$p$-coubure pour $p=1$.\\
 Apr\`es la courbure scalaire et la courbure sectionnelle usuelles, la courbure d'Einstein,
 est la plus importante parmi toutes les autres $p$-courbures.
Sa positivit\'e  a \'et\'e pr\'ec\'edemment \'etudi\'ee  mais seulement dans
le cadre
g\'en\'eral des autres $p$-courbures.
C'est d'ailleurs pour cette raison  qu'on lui a consacr\'e une \'etude \`a part 
dans \cite{Labbi6}.\\
La positivit\'e de la courbure d'Einstein  entra\^{i}ne celle de la courbure scalaire.
A l'heure actuelle, on connait  beaucoup de classes de vari\'et\'es admettant des m\'etriques \`a courbure
d'Einstein positive. En revanche, il nous manque  de trouver des exemples de vari\'et\'es
admettant des m\'etriques \`a courbure scalaire positive mais ne portant aucune m\'etrique 
\`a courbure d'Einstein positive. Les vari\'et\'es candidates sont $S^2\times T^n$, mais 
jusqu'\`a pr\'esent aucune d\'emonstration compl\`ete n'a pu \^etre \'etablie!\\
Rappelons   que de tels exemples doivent \^etre cherch\'es en dehors de la classe
 des vari\'et\'es compactes 2-connexe  et de dimension
$\geq 7$, voir \cite{Labbi3, Labbi6}.\\
Dans \cite{Labbi6}, on a 
 utilis\'e les techniques des surfaces minimales de Schoen-Yau dans l'ultime but de trouver
des obstructions topologiques \`a l'existence des m\'etriques \`a courbure d'Einstein positive, 
en dehors des
obstructions d\'ej\`a connues pour la courbure scalaire positive.\\
Dans deux articles c\'el\`ebres  \cite{Schoen-Yau} et \cite{SY2}, 
Schoen et Yau ont utilis\'e les
techniques des surfaces minimales pour obtenir des obstructions topologiques \`a l'\'existence
d'une  m\'etrique \`a courbure scalaire positive.
Le point cl\'e de leur d\'emarche est le r\'esultat suivant:
\begin{theorem}
 Soit $ (M,g)$ une vari\'et\'e compacte admettant une m\'etrique \`a courbure scalaire positive
 et telle que  $\dim M=m\geq 3$. Soit  $V$  une sous-vari\'et\'e compacte, lisse, de dimension
$(m-1)$, immerg\'ee dans $M$ et \`a fibr\'e normal trivial. Si de plus
$V$ est un minimum local de la
 $(m-1)$-volume, alors $V$ admet aussi une  m\'etrique \`a courbure scalaire positive.
 \end{theorem}
\noindent
Le r\'esultat suivant  est \'egalement post\'erieur \`a la  th\`ese o\`u on a
 g\'en\'eralis\'e le th\'eor\`eme pr\'ec\'edent au cas de la courbure d'Einstein comme suit:

\begin{theorem}[\cite{Labbi6}]
 Soit $ (M,g)$ une vari\'et\'e compacte admettant une m\'etrique \`a courbure 
d'Einstein positive
  et de dimension $m\geq 4$.
Soit  $V$  une sous-vari\'et\'e compacte, lisse, de dimension
$(m-2)$, immerg\'ee dans $M$ et \`a fibr\'e normal globalement plat. Si de plus
$V$ est un minimum local de la
 $(m-2)$-volume, alors $V$ admet  une m\'etrique \`a courbure scalaire positive.
\end{theorem}\noindent
Par  ``fibr\'e normal globalement plat'' on entend dire que le fibr\'e normal de $V$
poss\`ede deux sections globalement parall\`eles et partout orthonormales.
\par\medskip\noindent
Pour compl\'eter ce travail, il faudrait d\'emontrer
 l'existence de telles immersions pour certaines vari\'et\'es. En particulier, 
on ignore  toujours si on peut  immerger
le tore $T^2$ comme dans le th\'eor\`eme dans la vari\'et\'e produit $S^2\times T^2$. \\

 \section{Seconde courbure de Gauss-Bonnet-Weyl positive}
Dans toute la suite on \'ecrira en abr\'eg\'e {\bf  SCGBW} pour dire
{\bf seconde courbure de Gauss-Bonnet-Weyl}.
 Rappelons, comme on  vient de le voir dans le chapitre pr\'ec\'edent, que les vari\'et\'es d'Einstein ont leur {\sl SCGBW}
non-n\'egative et non identiquement nulle \`a moins qu'elles soient plates.\\
Par ailleurs il nous semble 
plausible,  comme dans le cas de la courbure scalaire, que toute vari\'et\'e portant une
m\'etrique riemannienne \`a {\sl SCGBW} non-n\'egative et non identiquement nulle admet aussi
une m\'etrique riemannienne \`a {\sl SCGBW} strictement positive.\\
Motiv\'e par ces faits,
on a \'etudi\'e dans \cite{Labbi8}, en vue d'une ultime classification, 
les vari\'et\'es admettant
une m\'etrique \`a {\sl SCGBW} positive.\\
Dans ce chapitre, on expose notre contribution \`a ce sujet.\\

\subsection{Liens avec la positivit\'e des autres courbures}
Le th\'eor\`eme suivant \'etabli une relation entre la positivit\'e de la  $p$-courbure
et celle de la  {\sl SCGBW}.:\par\medskip\noindent
\begin{theorem}[\cite{Labbi8}]
Soit $(M,g)$ une vari\'et\'e riemannienne de dimension  $n\geq 4$ et \`a $p$-courbure
non-n\`egative (resp. positive) telle que $p\geq {n\over 2}$.
 Alors, la  {\sl SCGBW} de $(M,g)$ est non-n\'egative (resp. positive).
 De plus, elle est identiquement nulle si et seulement si $(M,g)$ est plate.\end{theorem}
\par\medskip\noindent
Rappelons que la positivit\'e de la courbure sectionnelle  entra\^{\i}ne celle de la
 $p$-courbure pour $0\leq p\leq n-2$. De m\^eme que la courbure isotrope positive entra\^{\i}ne
 la positivit\'e de la $p$-courbure
pour
 $p\leq n-4$, voir
 \cite{Labbi4}.\\

Par cons\'equent, en dimensions sup\'erieures, la positivit\'e de la {\sl SCGBW} est plus
faible que celle de la courbure sectionnelle
ou de la courbure isotrope. Pr\'ecis\'ement, on a
\begin{corollary}[\cite{Labbi8}]
\begin{enumerate}
\item Si une vari\'et\'e riemannienne de dimension  $\geq 4$ est \`a courbure sectionnelle
non-n\'egative (resp.  positive) alors  sa {\sl SCGBW} est
non-n\'egative (resp.  positive).
 De plus, $h_4\equiv 0$
  si et seulement si elle est plate.
\item
Si une vari\'et\'e riemannienne de dimension  $\geq 8$ est \`a courbure isotrope
non-n\'egative (resp.  positive) alors sa  {\sl SCGBW} est
 non-n\'egative (resp.  positive).
 De plus, $h_4\equiv 0$
  si et seulement si elle est plate.
\end{enumerate}
\end{corollary}
Notons que la premi\`ere  partie du corollaire pr\'ec\'edent g\'en\'eralise un r\'esultat
de Milnor  en dimension 4.
\subsection{Cons\'equences}
Il d\'ecoule des  r\'esultats pr\'ec\'edents que les groupes de Lie munis d'une m\'etrique
 biinvariante et les vari\'et\'es riemanniennes normales homog\`enes ont leurs
{\sl SCGBW} non-n\'egative. De plus, en utilisant nos r\'esultats sur la $p$-courbure
 \cite{Labbi2,Labbi3} et le th\'eor\`eme pr\'ec\'edent on a montr\'e les
 cons\'equences suivantes:

\begin{corollary}[\cite{Labbi8}]
\begin{enumerate}
\item Soit $G$ un groupe de Lie  compact, connexe, de dimension $\geq 4$ et de rang  $r<[{\dim G+1\over 2}]$ muni
d'une m\'etrique biinvariante  $b$. Alors,  $(G,b)$ est \`a  {\sl SCGBW} positive.\par
\item  Si $G/H$ est une  vari\'et\'e riemannienne normale homog\`ene de dimension $\geq 4$
telle que le rang $r$ de $G$ satisfait
$r<[{\dim(G/H)+1\over 2}]$ alors elle est \`a {\sl SCGBW} positive.
\item
Si une vari\'et\'e compacte  $M$ de dimension $\geq 4$ admet une action lisse
d'un groupe de Lie  compact, connexe et de rang $r>[{\dim M+1\over 2}]$, alors 
 $M$ admet une m\'etrique riemannienne \`a {\sl SCGBW} positive.

\end{enumerate}
\end{corollary}

\subsection{Conjecture de Hopf alg\'ebrique g\'en\'eralis\'ee}
  Dans le cas o\`u la dimension  $n$ de la vari\'et\'e $M$ en question est paire,
la conjecture de Hopf alg\'ebrique stipule que la positivit\'e de la courbure sectionnelle
entra\^{\i}ne la positivit\'e
de l'int\'egrand de Gauss-Bonnet, i.e.  $h_{n}>0$.
Alors on peut se demander si plus g\'en\'eralement:
\par\smallskip
{\sl
La positivit\'e de la courbure sectionnelle
entra\^{\i}ne la positivit\'e
de $h_{2k}$, pour tout $k$ tel que  $2\leq 2k\leq n$?}\\
 Ou du moins dans le cas compact (mais \'eventuellement les deux probl\`emes sont 
\'equivalents!): \\
{\sl
La positivit\'e de la courbure sectionnelle
entra\^{\i}ne la positivit\'e
de la courbure totale $\int_M h_{2k}{\rm dvol}$, pour tout $k$ tel que  $2\leq 2k\leq n$?}\\
On sait maintenant que ceci est vrai pour  $k=1$ et $k=2$. Le probl\`eme reste ouvert pour
$k\geq 3$.\\
Il est opportun de mentionner que des contre-exemples purement alg\'ebriques
\`a cette conjecture existent, voir \cite{Geroch,Bourguignon2}. En revanche,
 on ignore  toujours s'il existe une v\'eritable
vari\'et\'e riemannienne \`a courbure sectionnelle positive sans que son int\'egrand de Gauss
Bonnet soit positive.

\subsection{Submersions riemanniennes et $h_4>0$}
On utilise ici la technique qui consiste \`a faire varier la m\'etrique de l'espace
total d'une submersion riemannienne pour faire appara\^ itre  la {\sl SCGBW} positive.
Signalons que cette technique a  auparavant \'et\'e utilis\'e par Cheeger, Lawson-Yau, Besse ...
Le r\'esultat obtenu est le suivant:
\begin{theorem}[\cite{Labbi8}]  Soit $M$ l'espace total   d'une submersion riemannienne.
Si $M$ est compact
et si les fibres (munis de la m\'etrique induite) sont \`a {\sl SCGBW} positive
alors la vari\'et\'e  $M$ admet une m\'etrique riemannienne \`a {\sl SCGBW} positive.
\end{theorem}
\par\medskip\noindent
Ce r\'esultat admet les cons\'equences directes suivantes:
\begin{corollary}
\begin{enumerate}
\item Le produit $S^p\times M$ d'une vari\'et\'e arbitraire $M$ compacte avec la
sph\`ere  $S^p,p\geq 4$ admet une m\'etrique riemannienne \`a {\sl SCGBW} positive.

\item Si une vari\'et\'e compacte  $M$ admet un feuilletage riemannien
tel que les feuilles sont \`a {\sl SCGBW} positive
alors la vari\'et\'e   admet une m\'etrique riemannienne \`a {\sl SCGBW} positive.
\item  Si une vari\'et\'e compacte  $M$  de dimension $\geq 4$ admet une action lisse et libre
d'un groupe de Lie $G$ compact, connexe et de rang $r$ tel que $r<[{\dim G+1\over 2}]$,
alors  $M$ admet une m\'etrique riemannienne \`a {\sl SCGBW} positive.
\end{enumerate}
\end{corollary}
\par\medskip\noindent
Il r\'esulte imm\'ediatement de la premi\`ere partie du corollaire pr\'ec\'edent qu'il
n'y a pas de restrictions sur le groupe fondamental d'une vari\'et\'e de dimension $\geq 8$
pour porter une m\'etrique riemannienne \`a {\sl SCGBW} positive. On verra ult\'erieurement
que
le m\^eme r\'esultat est vrai aussi pour les dimensions $\geq 6$.
\subsection{Chirurgies et $h_4>0$}
La stabilit\'e de la courbure scalaire positive sous chirurgies en codimensions $\geq 3$
a \'et\'e  le point cl\'e de la classification des vari\'et\'es simplement connexes
 \`a courbure scalaire
 positive d\^u
\`a Gromov-Lawson \cite{Gromov-Lawson} et Schoen-Yau \cite{Schoen-Yau} et Stolz \cite{Stolz}.
On a montr\'e dans \cite{Labbi8} un r\'esultat similaire pour la  {\sl SCGBW}:
\begin{theorem}[\cite{Labbi8}] Si une vari\'et\'e  $M$ est obtenue \`a partir d'une vari\'et\'e compacte
$X$ par chirurgies en codimension
 $\geq 5$, et si $X$ admet  une m\'etrique riemannienne \`a {\sl SCGBW} positive,
alors $M$  admet aussi une m\'etrique riemannienne \`a {\sl SCGBW} positive.\\
En particulier, la somme connexe de deux vari\'et\'es compactes de dimensions $\geq 5$
 portant chacune
une  m\'etrique riemannienne \`a {\sl SCGBW} positive,
 porte aussi une m\'etrique \`a {\sl SCGBW} positive.
\end{theorem}
Comme cons\'equence de ce  th\'eor\`eme, on a montr\'e qu'il n'y a pas
de restrictions sur le groupe fondamental d'une vari\'et\'e de dimension $\geq 6$
 \`a {\sl SCGBW} positive. Pr\'ecis\'ement, on a
\begin{corollary}[\cite{Labbi8}] Soit  $G$ un groupe \`a pr\'esentation finie.
Alors pour tout  $n\geq 6$, il existe une vari\'et\'e compacte de dimension $n$
portant une m\'etrique \`a {\sl SCGBW} positive et telle que  $\pi_1(M)=G$.
\end{corollary}

\section{Courbure isotrope positive}
\subsection{Introduction}
La courbure isotrope a \'et\'e introduite par Micallef et Moore \cite{Moore}  pour les 
vari\'et\'es
de dimension $\geq 4$. Cette courbure joue le m\^eme r\^ole dans l'\'etude de
la stabilit\'e des 2-sph\`eres
harmoniques et  des surfaces minimales, que celui exerc\'ee par
 la courbure sectionnelle
dans l'\'etude de la stabilit\'e des g\'eod\'esiques.\\
L'existence d'une m\'etrique riemannienne \`a courbure isotrope positive
sur une vari\'et\'e compacte entra\^{\i}ne l'annulation des groupes d'homotopie
$\pi_i(M)$ pour $2\leq i\leq [n/2]$, voir \cite{Moore}. En particulier, si la vari\'et\'e est en plus
simplement connexe alors elle doit \^etre hom\'eomorphe \`a une sph\`ere.\\
Il y a une similarit\'e entre la courbure isotrope et la $(n-4)$-courbure, o\`u $n$ est la dimension
de la vari\'et\'e en question.
En effet, elles se distinguent seulement par un terme en la courbure de Weyl, et par cons\'equent, 
elles co\"{\i}ncident dans le cas conform\'ement plat. Cette analogie entre ces deux courbures est 
le point de d\'epart de notre contribution  \cite{Labbi4,Labbi5}
dans l'\'etude de  positivit\'e de la courbure isotrope. Nous allons dans un premier temps rappeler
la d\'efinition de la courbure isotrope.\\
Soit $(M,g)$ une vari\'et\'e riemannienne de dimension $n$.
Pour tout  $m\in M$, le produit scalaire  $g$
 sur l'espace tangent  $T_mM$,
peut  s'\'etendre de deux mani\`eres au complexifi\'e
$T_mM\otimes {\bf C}$:\par
- En tant qu'une forme bilin\'eaire complexe, qu'on notera par  $g(.,.)$.\par

- En tant qu'un produit scalaire Hermitien, qu'on notera par
 $<.,.>$.\par\medskip\noindent

Soit $ R:\wedge^2 M\rightarrow \wedge^2 M$ l'op\'erateur de courbure de
$(M,g)$. On notera aussi par $R$ son extension lin\'eaire complexe \`a
 $\wedge ^2M \otimes
{\bf C}$.\par\smallskip\noindent
On associe \`a  tout 2-plan complexe  $P\subset T_mM \otimes
{\bf C}$  sa courbure sectionnelle complexe  $K_C(P)$
d\'efinie par
 $$K_C(P)={< R(v\wedge w),v\wedge w> \over ||v\wedge w||^2},$$
o\`u $\{v,w\}$ est une base quelconque de  $P$.
On dit qu'un sous-espace vectoriel complexe  $P\subset T_mM \otimes
{\bf C}$ est isotrope si  $g(v,v)=0$ pour tous les vecteurs  $v\in P$.
 \par\medskip\noindent
{\sc Definition.} On dit qu'une vari\'et\'e riemannienne  $(M,g)$ est \`a {\sl courbure
isotrope positive} au pont $m\in M$ si
$K_C(P)>0$ pour tout les 2-plans complexes et isotropes dans
$T_mM\otimes {\bf C}$.\par\medskip\noindent
Cette d\'efinition est \'equivalente \`a la condition suivante sur le tenseur de courbure:
\begin{equation}
K(e_1,e_3)+K(e_1,e_4)+K(e_2,e_3)
+K(e_2,e_4)>2 R(e_1,e_2,e_3,e_4)
\end{equation}
Et ceci pour tout les vecteurs  orthonorm\'es $\{e_1,e_2,e_3,e_4\}$
dans l'espace tangent  $T_mM$.\\
Une propri\'et\'e importante de la courbure isotrope positive est qu'elle est conserv\'e
par le flot de Ricci. En effet, Hamilton \cite{Hamilton} a utilis\'e le flot de Ricci pour classifier
les vari\'et\'es compactes de dimension 4 admettant une m\'etrique \`a courbure isotrope
positive et n'ayant pas d'``espaces formes essentiels et incompressibles". Pr\'ecis\'ement
il a d\'emontr\'e que de telles vari\'et\'es sont diff\'eomorphes soit \`a
$S^4$, $RP^4$, $S^1\times S^3$, $S^1\tilde{\times} S^3$ soit \`a une somme connexe des vari\'et\'es pr\'ec\'edentes.

\subsection{Obstructions}
L'existence d'une m\'etrique riemannienne \`a courbure isotrope positive
sur une vari\'et\'e compacte entra\^{\i}ne des restrictions importantes sur sa
 topologie.
En effet,  Micallef
et Moore \cite{Moore} ont d\'emontr\'e que:
\par\medskip
 {\sl  Soit  $M$ une vari\'et\'e riemannienne compacte de
dimension
$n\geq 4$. Si $M$ est \`a courbure isotrope positive alors
$\pi_i(M)=\{0\}$ pour $2\leq i\leq [n/2]$. En particulier, si
$M$ est de plus simplement
connexe, alors
$M$ est hom\'eomorphe \`a une sph\`ere.}\par\medskip\noindent
D'autre part,  Micallef et Wang \cite{Micallef-Wang},
et Seaman \cite{Seaman} ont d\'emontr\'e l'annulation du second nombre de Betti dans le cas o\`u
la dimension de la vari\'et\'e est paire. On croit que tous les nombres de Betti
$b_k:2\leq k\leq n-2$
devraient s'annuler, mais pour le moment on n'a pas de d\'emonstration de cette conjecture.
  En revanche, on sait d\'emontrer ce fait  dans le cas conform\'ement plat
 \cite{Labbi4}. Le r\'esultat
est le suivant:
\begin{theorem}[\cite{Labbi4}]
\begin{enumerate}
\item Soit   $(M,g)$ une vari\'et\'e  compacte conform\'ement plate de dimension $n$  et portant une m\'etrique \`a courbure isotrope positive
alors  $H^m(M,{\bf R}) = 0$
pour
 $ 2 \leq  m\leq n-2  $.
\item Soit $(M,g)$ une vari\'et\'e  compacte conform\'ement plate de dimension $n$
\`a courbure isotrope  non-n\'egative  et telle que
   $H^2(M,{\bf R}) \not= 0$ alors soit $(M,g)$ est plate soit elle est revetue par
$S^{n-2}\times H^2$.
\end{enumerate}
\end{theorem}
Notons que des  r\'esultats du m\^eme type ont  aussi \'et\'e  d\'emontr\'es par  Nayatani
 \cite{Nayatani}
et Mercuri-Noronha \cite{Noronha}.\\
Finalement,  l'exemple de $S^1\times S^n$, $n\geq 3$, montre que le groupe fondamental d'une vari\'et\'e
\`a courbure isotrope positive peut \^etre infini. En plus, parce que la courbure isotrope
positive est conserv\'ee par l'op\'eration de sommes connexes, le groupe fondamental peut
\^etre assez grand. En contre partie, Fraser \cite{Fraser} a d\'emontr\'e, en utlisant les techniques
de surfaces minimales, la r\'estriction suivante sur le groupe fondamental:
\\
 {\sl Soit $M$ une vari\'et\'e riemannienne compacte de dimension $n\geq 5$
\`a courbure isotrope positive. Alors, le groupe fondamental de $M$ ne contient pas de 
sous-groupes isomorphes \`a $\mathbb{Z}\oplus \mathbb{Z}$}.\\

\subsection{Constructions}

Par ailleurs, concernant les  constructions de m\'etriques \`a courbure isotrope positive,
  Micallef et Wang \cite{Micallef-Wang}  ont
 d\'emontr\'e le r\'esultat suivant:
\\
{\sl La somme connexe de deux vari\'et\'es  de dimensions  $\geq 4$ et
portant chacune  une m\'etrique \`a courbure isotrope positive, porte aussi une
m\'etrique \`a courbure isotrope positive.} \\
En utilisant les techniques de submersions riemanniennes,
on a  d\'emontr\'e les constructions
suivantes:
\begin{theorem}[\cite{Labbi4}]
Soit
 $\pi : (M,g)\rightarrow S^1$ une submersion riemannienne  sur le cercle, o\`u  $M$
 d\'esigne une vari\'et\'e riemannienne compacte de dimension $\geq
 4$.
\\
 Supposons que les fibres de 
$\pi$ (munis de la m\'etrique induite) satisfont la condition  de positivit\'e (A) ci-dessous,
alors,  $M$ porte une m\'etrique \`a courbure isotrope positive.
\end{theorem}
O\`u  (A) d\'esigne la condition suivante sur le tenseur de courbure:
\par\medskip
$$  K(e_j,e_k)+K(e_j,e_l) >\mid
R(e_i,e_j,e_k,e_l)\mid, \eqno {(A)}$$
pour tous les vecteurs tangents orthonorm\'es  $\{e_i,e_j,e_k,e_l\}$.
\par\medskip\noindent
Notons que, la c\'el\`ebre condition  de ${1\over 4}$-pincement stricte qui porte sur la courbure
sectionnelle
  entra\^{\i}ne la condition de positivit\'e  (A).  De plus, la  condition (A) entra\^{\i}ne 
en m\^eme temps
la positivit\'e de la courbure isotrope et celui de la courbure de Ricci.
Enfin, en dimension 3, la condition (A) est \'equivalente \`a la positivit\'e de la
 courbure de Ricci.
En particulier, on obtient:
\begin{corollary}[\cite{Labbi4}]
 Soit  $\pi : (M,g)\rightarrow S^1$ une submersion riemannienne 
sur le cercle, o\`u  $M$
 d\'esigne une vari\'et\'e riemannienne compacte de dimension $4$.\\
 Supposons que les fibres de
$\pi$ (munis de la m\'etrique induite) sont \`a courbure de Ricci positive
alors $M$ porte une m\'etrique \`a courbure isotrope positive.\end{corollary}
Comme cons\'equence des r\'esultats pr\'ec\'edents on obtient les exemples suivants
de vari\'et\'es admettant des m\'etriques \`a courbure isotrope positive:\\
\begin{enumerate}
\item Soit  $F$ une vari\'et\'e de dimension  $\geq 3$ portant une m\'etrique riemannienne
$g$ qui satisfait la condition de positivit\'e (A).
Soit $\phi \in Isom (F,g)$ et soit
$\rho : {\bf{Z}} \longrightarrow Isom(F\times {\bf R})$
d\'efinie par
$$n\longrightarrow \phi_n(x,t)=(\phi^n(x),t+n).$$
La vari\'et\'e  $M={{F\times {\bf R}}\over \rho}$
est alors l'espace total d'une submersion riemannienne et satisfait
les hypoth\`eses du th\'eor\`eme.
Elle admet donc une m\'etrique \`a courbure isotrope positive.
\item Soit $M$ une vari\'et\'e compacte  de dimension $\geq 4$ et admettant un feuilletage
riemannien
de  codimension 1
 tel que les feuilles, munis de la m\'etrique induite,
satisfont la condition
 (A), alors $M$ admet  une m\'etrique \`a courbure isotrope positive.
\item Si une vari\'et\'e compacte de dimension 4 admet une action libre et lisse du groupe
$SU(2)$ ou $SO(3)$ alors elle admet  une m\'etrique \`a courbure isotrope positive.
\end{enumerate}
Remarquons  que le dernier r\'esultat  ne se g\'en\'eralise pas au
 cas d'une action non libre.
En effet,  $S^2\times S^2$ admet une action effective du groupe  $SO(3)$, mais elle n'admet
aucune m\'etrique  \`a courbure isotrope positive.\\
\subsection{Chirurgies et courbure isotrope positive}
Du fait de l'analogie existante entre courbure isotrope et $(n-4)$-courbure, on peut
esp\'erer que les r\'esultats de chirurgies d\'emontr\'es pour la $p$-courbure \cite{Labbi3}
se g\'en\'eralisent au cas de la courbure isotrope positive.
Pr\'ecis\'ement, on peut s'attendre \`a ce que
la courbure isotrope positive soit stable sous chirurgies en codimension $\geq n-1$.
D'autant plus qu'une telle stabilit\'e a \'et\'e annonc\'ee, sans d\'emonstration, 
dans \cite{Gromov2}.
 Cette question a fait l'objet de
l'article \cite{Labbi5}. Dans lequel on a d\'emontr\'e la stabilit\'e de la courbure isotrope positive
sous une chirurgie de codimension $n$, i.e. sous sommes connexes. Ce r\'esultat a \'et\'e premi\`erement d\'emontr\'e
par Micallef-Wang en utilisant les techniques de  Schoen-Yau \cite{Schoen-Yau}. Notre d\'emonstration
\'etait \`a la mani\`ere de Gromov-Lawson \cite{Gromov-Lawson}.\\
Mais  malheureusement, la stabilit\'e sous une chirurgie en codimension $n-1$ n'a pas
en g\'en\'eral toujours lieu.
En effet, on montre qu'un tel r\'esultat entra\^{\i}nerai alors que tout groupe
\`a pr\'esentation finie
peut \^etre r\'ealis\'e comme le groupe fondamental d'une vari\'et\'e \`a courbure isotrope
 positive (de dimension suffisamment grande). Ce qui vient
contredire le r\'esultat de Fraser cit\'e ci-dessus.
\bigskip\noindent

\newpage 
\part{Perspectives de recherches}
\vspace{1cm}
\begin{center}
{\sl
Comme dans tout travail de recherche math\'ematiques, le point de d\'epart est un probl\`eme
soulev\'e au bout duquel se ramifient d'autres questions naturelles.
L'arbre porte un r\^eve au bout de chaque branche et le chercheur avance
au gr\'e des rameaux. Si l'on ignore le lieu et le temps o\`u on atteindra le fruit, on  a
au moins la certitude de notre apport en oxyg\`ene \`a l'arbre ramifi\'e de la recherche
et aussi d'en \^etre un \'el\'ement.}
\end{center}
\vspace{0.5cm}
Les ramifications naturelles \`a mes recherches sont les suivantes:\\

\section*{A. Obstructions topologiques pour $h_4>0$:}
Un des probl\`emes de grande importance en g\'eom\'etrie contemporaine est l'\'etude des vari\'et\'es
d'Einstein. Ce projet se situe dans cette direction puisque les vari\'et\'es d'Einstein
ont leur  seconde courbure de Gauss-Bonnet-Weyl non-n\'egative et non identiquement nulle 
\`a moins
qu'elles soient plates. Et ceci, rappelons le, ind\'ependement du signe de la constante
 d'Einstein.\\
Dans ce projet, il s'agit d'abord de d\'emontrer  qu'une vari\'et\'e compacte admettant une m\'etrique
\`a courbure $h_4\geq 0$ non identiquement nulle porte aussi une m\'etrique \`a $h_4>0$.\\
Ensuite, on doit chercher des obstructions topologiques 
(comme dans le cas de la courbure scalaire) \`a l'existence des m\'etriques \`a $h_4>0$.\\
De tels r\'esultats seraient  tr\`es int\'eressants dans le sens o\`u ils prouveraient l'existence
de vari\'et\'es en dimensions sup\'erieures sans m\'etriques d'Einstein.\\
Une autre voie possible de  recherche pour ces obstructions, est de consid\'erer l'analogue 
poly\`edral de $h_4$. C'est une mesure port\'ee par le squelette de codimension $4$ 
 du poly\`edre en question. Pour des approximations poly\`edrales convenables d'une vari\'et\'e
riemannienne,
ces mesures discr\`etes convergent en mesure vers l'int\'egrand $h_4{\rm dvol}$. Voir les
 travaux de
Cheeger, M\"{u}ller et Schr\"{a}der \cite{C-M-S} et le compte-rendu de Lafontaine \cite{Lafontaine-poly}.\\
Il s'agit donc d'\'etudier les propri\'et\'es de positivit\'e de ces analogues poly\`edraux. En 
particulier, de
voir si on peut  avoir des th\'eor\`emes d'annulation des nombres de Betti.

\section*{B. Le probl\`eme de  Yamabe g\'en\'eralis\'e:} 
Ce probl\`eme  a d\'ej\`a \'et\'e  abord\'e dans ce texte. Le lecteur de ces lignes est invit\'e
\`a se  r\'ef\`erer donc au chapitre 2 pour les
motivations. Il s'agit ici d'\'etudier  
la conjecture suivante:\\
{\sl Dans  toute classe conforme d'une m\'etrique
riemannienne sur une vari\'et\'e $C^\infty $ compacte donn\'ee, il existe une m\'etrique 
riemannienne
\`a courbure de
Gauss-Bonnet-Weyl   $h_{2k}$ constante.}\\
\section*{C. Th\'eor\`eme d'annulation pour la courbure isotrope positive}
 L'objectif  ici est de d\'emontrer que la positivit\'e de la courbure isotrope sur une vari\'et\'e
compacte de dimension $n\geq 4$
  entra\^{\i}ne l'annulation des nombres
 de Betti $b_i$ de la vari\'et\'e  pour $2\leq i\leq n-2$.\\
Ce probl\`eme est motiv\'e d'une part par l'analogie existante entre courbure isotrope
et $(n-4)$-courbure et d'autre part par le fait que le tenseur de $p$-courbure appara\^{\i}t
 naturellement
dans l'op\'erateur de courbure de Weitzenb\"{o}ck. Signalons qu'on a d\'ej\`a d\'emontr\'e ce 
 r\'esultat pour une m\'etrique conform\'ement
plate. 

\section*{D. Conjecture de Hopf alg\'ebrique g\'en\'eralis\'ee}
Ce probl\`eme  a d\'ej\`a \'et\'e  abord\'e et motiv\'e ant\'erieurement dans ce projet. 
Je me contenterai ici  du simple \'enonc\'e de cette perspective de recherche:  \\
{\sl
La positivit\'e de la courbure sectionnelle
entra\^{\i}ne t-elle la positivit\'e
de $h_{2k}$, pour tout $k$ tel que  $2\leq 2k\leq n$? o\`u $n$ d\'esigne la dimension de la 
vari\'et\'e en question.}\\
 Dans le cas compact, on esp\`ere d\'emontr\'e  une version faible de ce probl\`eme
 (mais \'eventuellement les deux probl\`emes sont 
\'equivalents!): \\
{\sl
Pour une vari\'et\'e compacte, la positivit\'e de la courbure sectionnelle
entra\^{\i}ne la positivit\'e
de la courbure totale $\int_M h_{2k}{\rm dvol}$, pour tout $k$ tel que  $2\leq 2k\leq n$?}\\
Ceci \'etant  vrai pour  $k=1$ et $k=2$, voir chapitre 2. Le probl\`eme reste ouvert pour
$k\geq 3$. 

\section*{E.  Vari\'et\'es $(p,q)$-Einstein:} Le but de cette recherche serait d'\'etudier 
les g\'en\'eralisations naturelles suivantes
des m\'etriques d'Einstein.\\
Pour $1\leq p<2q$, on dira qu'une vari\'et\'e riemannienne de tenseur de courbure $R$
  est $(p,q)$-Einstein
 si son tenseur de Ricci g\'en\'eralis\'e  $c^{2q-p}R^q$ est proportionnel \`a la m\'etrique 
$g^p$. On retrouve les 
vari\'et\'es d'Einstein usuels pour $p=q=1$ et les vari\'et\'es $(2q)$-Einstein (discut\'ees 
dans le chapitre 2) pour $p=1$.\\
On connait d\'ej\`a une caract\'erisation de cette condition:
\begin{proposition}[\cite{Labbi7}] Pour $1\leq p<2q$, une m\'etrique riemannienne est 
$(p,q)$-Einstein si et seulement si pour tout $1\leq i\leq \min\{p,n-p\}$ on a
$\omega_i=0$. O\`u les $\omega_i$ sont les composantes de $R^q$ suivant la d\'ecomposition
$R^q=\sum_{i=0}^{2q}g^{2q-i}\omega_i$.
\end{proposition}
Pour $p=q=1$, on retrouve la caract\'erisation classique des vari\'et\'es d'Einstein usuelles
par l'annulation de la composante $\omega_1$ du tenseur de courbure de Riemann.\\
Aussi, il en  r\'esulte imm\'ediatement de notre \'etude  que si une vari\'et\'e compacte 
de dimension $n=4q$
est $(2q-1,q)$-Einstein alors sa caract\'eristique d'Euler-Poincarr\'e est non-n\'egative.\\
Ce dernier r\'esultat g\'en\'eralise une propri\'et\'e bien connue pour $q=1$.\\
Remarquons que la condition $(p,q)$-Einstein entraine $(p-1,q)$-Einstein pour tout $p\geq 2$.\\
Pa ailleurs, si la dimension de la vari\'et\'e est telle que $n=2q$, 
alors toute vari\'et\'e est $(p,q)$-Einstein. Donc, une telle
condition n'est pas vide que si $2q<n$.\\
Une question qui me semble int\'eressante serait de trouver une classification 
des vari\'et\'es qui sont
$(p,q)$-Einstein pour tout $1\leq q<n/2$ et tout  $1\leq p<2q$ (il suffit de prendre $p=2q-1$).\\
Cette question g\'en\'eralise, aux dimensions sup\'erieures, le probl\`eme
 de classification des vari\'et\'es  d'Einstein
en dimension 4.\\
Notons finalement  que les m\'etriques en question ici ne sont pas, en g\'en\'eral,
 critiques au sens
de Bleecker \cite{Bleeker}.

\section*{F. Vari\'et\'es conform\'ement plates g\'en\'eralis\'ees:} Comme dans la perspective
de recherche pr\'ec\'edente, l\`a aussi
les questions ne sont pas tr\`es pr\'ecises. Il s'agit d'\'etudier une classe de m\'etriques 
g\'en\'eralisant les m\'etriques conform\'ement plates.\\
Rappelons qu'une m\'etrique est conform\'ement plate (en dimension $\geq 4$)
 si son tenseur de courbure de Riemann
est divisible par la m\'etrique. Elles sont caract\'eris\'ees par l'annulation de la composante
$\omega_2$ dans la d\'ecomposition orthognale en composantes irr\'eductibles de $R$.
 En particulier, une vari\'et\'e  d'Einstein 
 conform\'ement plate est n\'ecessairement \`a courbure constante. \\
Par analogie, on dira qu'une vari\'et\'e riemannienne  de dimension $n$ 
est  de classe  $C(p,q)$ si son tenseur
de Gauss-Kronecker $R^q$ est divisible par $g^p$. Pour $p=q=1$, on 
retrouve les m\'etriques
conform\'ement plates usuelles. Pour $p=1$ et $q\geq 1$, cette condition est invariante sous 
un changement conforme de la m\'etrique. De telles m\'etriques sont  alors dites
$q$-conform\'ement plates et ont \'et\'e partiellement \'etudi\'ees par Nasu dans 
\cite{Nasu2,Nasu-Kojima}. \\
Les m\'etriques de classe $C(p,q)$ sont caract\'eris\'ees par l'annulation des  composantes
$\omega_i$ dans la d\'ecomposition orthogonale en composantes irr\'eductibles de
 $R^q=\sum_{i=0}^{2q}g^{2q-i}\omega_i$ pour 
$2q-p< i\leq 2q$.\\
On d\'eduit alors que, pour tout $q$ telle que $2q<n$, si une m\'etrique
 est \`a la fois $(q,q)$-Einstein (voir la section ci-dessus) et de classe  $C(q,q)$ alors 
elle est \`a $q$-courbure 
sectionnelle constante. \\
\section*{G. Obstructions pour la courbure d'Einstein positive}
La question ici est de trouver des classes de vari\'et\'es portant des
m\'etriques \`a courbure scalaire positive mais qui n'admettent aucune
m\'etrique \`a courbure d'Einstein positive.
La vari\'et\'e produit $S^2\times T^2$,   et plus g\'en\'eralement le
produit $S^2\times T^n$, $n\geq 2$ sont des candidats naturels.

\section*{H. Multiformes}
Cette perspective de recherche se pr\'esente comme la suite de mes travaux 
sur les formes doubles.\\
Rappelons q'une  forme double  $\omega\in \Lambda^{*p}V \otimes \Lambda^{*q}V$ peut \^etre 
consid\'er\'ee
comme une $p$-forme ordinaire  mais
\`a valeurs vectorielles dans $\Lambda^{*q}V$. Cette consid\'eration
 permet d'\'etudier les doubles formes \`a partir des formes. Elle
 a \'et\'e utilis\'ee par de nombreux auteurs, je cite par exemple
\cite{Besse,Bourguignon} .\\
Par ailleurs, certains de nos r\'esultats sur les  formes doubles peuvent
 \^etre g\'en\'eralis\'es au cas des  formes doubles \`a valeurs vectorielles.
 Un exemple de telles
formes doubles est
la seconde forme fondamentale d'une sous-vari\'et\'e de codimension $\geq 2$.\\
Ceci  permettrait d'\'etudier les "formes triples". Ces derniers \'etant les \'el\'ements
de l'espace $\Lambda^{*p}V \otimes \Lambda^{*q}V \otimes \Lambda^{*r}V$.\\
Il serait alors int\'er\'essant d'\'etudier en toute g\'en\'eralit\'e les
multiformes i.e. les \'el\'ements des
produits $\Lambda^{*{p_1}}V \otimes \Lambda^{*{p_2}}V \otimes ...\otimes \Lambda^{*{p_r}}V$.
En particulier, l'op\'erateur de Hodge peut \^etre g\'en\'eralis\'e de  $2^r-1$ diff\'erentes 
mani\`eres aux
multiformes d'ordre $r$ (en le faisant op\'erer \`a chaque fois sur diff\'erents 
facteurs du produit tensoriel).\\
Notons que tout tenseur de rang arbitraire peut \^etre vu comme une multiforme.
Certains aspects de ces probl\`emes ont \'et\`e r\'ecemment \'etudi\'es par
le physicien th\'eoricien J. M. M. Senovilla 
 dans une s\'erie
d'articles.
 En particulier, il a utlis\'e des extensions du produit de
Hodge et du produit int\'erieur aux multiformes pour d\'efinir, \`a partir d'un tenseur $T$
donn\'e,
des tenseurs sym\'etriques dites
de super-\'energie. Ces tenseurs sont quadratiques en $T$ et
 g\'en\'eralisent  les tenseurs dits de Bel et de Bel-Robinson, voir
\cite{Senovilla} et les ref\'erences qui y sont cit\'ees.

\section*{I. Courbures n\'egatives}
Suite aux travaux de Lohkamp \cite{Lohkamp}, on sait maintenant que toute vari\'et\'e
compacte admet des m\'etriques \`a courbure de Ricci n\'egative et donc en particulier 
des m\'etriques
\`a courbure scalaire n\'egative. Il serait alors int\'eressant de prouver la m\^eme
propri\'et\'e pour la seconde courbure de Gauss-Bonnet-Weyl, courbure d'Einstein, $p$-courbures,
 courbure isotrope et pour les courbures sectionnelles de Weitzenb\"{o}ck.

\newpage
\addcontentsline{toc}{part}{Appendix}

\appendix
\section{Annexe: Remarques sur la g\'eom\'etrie des tubes}
Durant une op\'eration de chirurgie sur une vari\'et\'e riemannienne
lelong d'une sous-vari\'et\'e
plong\'ee \cite{{Gromov-Lawson},{Labbi3},{Labbi8}}, on remplace la m\'etrique (ancienne)
induite par l'application  exponentielle  sur un voisinage  tubulaire de la
sous-vari\'et\'e plong\'ee
 par la m\'etrique de  Sasaki. Dans cet  appendice  on pr\'esente les r\'esultats de
\cite{Labbi10} o\`u on a \'etudi\'e  le
 comportement de ces deux m\'etriques au voisinage de la section nulle du
 fibr\'e normal.
\\
Soit  $(X,g)$ une vari\'et\'e $C^\infty$ riemannienne de  dimension $n+p$ et
$M$ une sous-vari\'et\'e (compacte)
plong\'ee
de dimesion $n$ dans $X$. Soit
 $$ T_\epsilon=\{(x,v): x\in M,v\in N_xM \quad {\rm and}\quad
g(v,v) <\epsilon^2\}$$
un tube de rayon  $\epsilon$ autour de  M.
Il existe $\epsilon_0>0$ tel que
${\rm exp}:T_\epsilon\rightarrow X$, est un diff\'eomorphisme sur son image
  pour tout $\epsilon\leq \epsilon_0$. 
Soit ${\rm exp}^*g$ le
  pull-back  de  la m\'etrique $g$ sur $X$.
\par
Le sous-fibr\'e  normal  $T_\epsilon$ peut \^etre aussi muni d'une deuxi\`eme
 m\'etrique naturelle, \`a savoir celle  de Sasaki.  C'est la m\'etrique
 $h$ compatible avec la connexion normale
telle que la projection naturelle $ \pi :(T_\epsilon,h)\rightarrow (M,g)$ soit une
  submersion riemannienne.\par
  \medskip
Soit $(p,rn)\in T_\epsilon$, o\`u $r<\epsilon_0$
 et $n$  un vecteur normal unitaire en $p$.
  Notons par $A_n$  l'application de Weingarten  de la sous-vari\'et\'e
    $M$  d\'efinie dans la
direction de  $n$.
  \par\medskip\noindent
Le r\'esultat suivant montre que ces deux m\'etriques sont en g\'en\'erale tangentes que
 jusqu'\`a
l'ordre 1. Elles le sont \`a l'ordre deux si et seulement si la sous-vari\'et\'e $M$ est
totalement g\'eod\'esique dans $(X,g)$.\\
\begin{theorem}[\cite{Labbi10}] Soit   $R$ le tenseur de courbure de Riemann
 de $(X,g)$. Alors pour tout $u_1,u_2\in T_{(p,rn)}T_\epsilon $, on a
\begin{equation*}
\begin{split}
{\rm exp}^*g(u_1,u_2)&=h(u_1,u_2)-2g(A_n\pi_*u_1,\pi_*u_2)r\\
+&\bigl\{g(A_n\pi_*u_1,A_n\pi_*u_2)+R(\pi_*u_1,n,\pi_*u_2,n)+
{2\over 3}R(\pi_*u_1,n,Ku_2,n)\\
+&{2\over 3}R(\pi_*u_2,n,Ku_1,n)+
{1\over 3}R(Ku_1,n,Ku_2,n)\bigr\}r^2+O(r^3).
\end{split}
\end{equation*}
O\`u $\pi_*$ et $K$ repr\'esentent r\'esp\'ectivement la diff\'erentielle de $\pi$ et 
l'application de connexion.
\end{theorem}
 \par\medskip\noindent
{\sc Remarque.} 
Notons que dans  \cite{Gromov-Lawson}, au d\'ebut de la preuve du lemme 2
page 430, il a \'et\'e affirm\'e que les m\'etriques  ${\rm exp}^*g$ et $h$
sont suffisamment proches dans la  topologie $C^2$ quand
 $r\rightarrow 0$. La m\^eme erreur est dans  \cite{Labbi3}. Ceci n'affecte pas
les r\'esultats correspondants
dans les deux articles, voir
 \cite{Labbi8}. \par\noindent
Une mani\`ere rapide pour ce rendre compte de cette propri\'et\'e est comme suit:
pour la m\'etrique $h$, la  section nulle  $M \hookrightarrow T_{\epsilon}$ est
totalement g\'eod\'esique (car pour une submersion riemannienne  le relev\'e horisontal
d'une g\'eod\'esique est une g\'eod\'esique).
 D'autre part, la  section nulle  $M \hookrightarrow T_{\epsilon}$ est
totalement g\'eod\'esique
pour la m\'etrique ${\rm exp}^*g$ si et seulement si  $M$ est totalement g\'eod\'esique
dans  $(X,g)$.\\
\par\smallskip\noindent
Dans le cas o\`u  l'espace ambiant  $(X,g)$ est l'espace  Euclidien
${\mathbb R}^n$,  on montre la relation  simple suivante entre les deux m\'etriques:

  \begin{theorem}[\cite{Labbi10}]
 Soit  $M$ une sous-vari\'et\'e plong\'ee dans l'espace Euclidien,
alors pour tout  $u_1,u_2\in T_{(p,rn)}T_\epsilon $, on a
\begin{equation*}
{\rm exp}^*g(u_1,u_2)=h(u_1,u_2)-2g(A_n\pi_*u_1,\pi_*u_2)r
+g(A_n\pi_*u_1,A_n\pi_*u_2)r^2.
\end{equation*}
\end{theorem}
On  d\'emontre des r\'esultas similaires dans le cas o\`u  $(X,g)$ est \`a courbure
constante arbitraire, voir \cite{Labbi10}.

\vspace{2cm}
\noindent
Labbi Mohammed-Larbi\\
  Department of Mathematics,\\
 College of Science, University of Bahrain,\\
  P. O. Box 32038 Bahrain.\\
  E-mail: labbi@sci.uob.bh

\end{document}